\documentclass[12pt,english]{amsart}
\usepackage[T1]{fontenc}
\usepackage[latin9]{inputenc}
\usepackage{geometry}
\usepackage{verbatim}
\usepackage{amstext}
\usepackage{amsthm}
\usepackage{amssymb}

\makeatletter
\numberwithin{equation}{section}
\numberwithin{figure}{section}
\theoremstyle{plain}
\newtheorem{thm}{\protect\theoremname}[section]
\theoremstyle{definition}
\newtheorem{defn}[thm]{\protect\definitionname}
\theoremstyle{remark}
\newtheorem{rem}[thm]{\protect\remarkname}
\theoremstyle{plain}
\newtheorem{cor}[thm]{\protect\corollaryname}
\theoremstyle{plain}
\newtheorem{lem}[thm]{\protect\lemmaname}
\theoremstyle{plain}
\newtheorem{fact}[thm]{\protect\factname}

\usepackage{babel}
\usepackage[unicode=true]
 {hyperref}
\usepackage{breakurl}

\makeatother

\usepackage{babel}
\providecommand{\corollaryname}{Corollary}
\providecommand{\definitionname}{Definition}
\providecommand{\factname}{Fact}
\providecommand{\lemmaname}{Lemma}
\providecommand{\remarkname}{Remark}
\providecommand{\theoremname}{Theorem}

\begin{document}
\title[Constants in the John-Nirenberg inequality]{Unions of cubes in $\mathbb{R}^{n}$, combinatorics in $\mathbb{Z}^{n}$
and the John-Nirenberg and John-Str\"omberg inequalities}
\author{Michael Cwikel}
\address{Department of Mathematics, Technion - Israel Institute of Technology,
Haifa 32000, Israel }
\email{mcwikel@math.technion.ac.il}
\thanks{This research was supported by funding from the Martin and Sima Jelin
Chair in Mathematics, by the Technion V.P.R.~Fund and by the Fund
for Promotion of Research at the Technion.}
\begin{abstract}
Suppose that the $d$-dimensional unit cube $Q$ is the union of three
disjoint ``simple'' sets $E$, $F$ and $G$ and that the volumes
of $E$ and $F$ are both greater than half the volume of $G$. Does
this imply that, for some cube $W$ contained in $Q$. the volumes
of $E\cap W$ and $F\cap W$ both exceed $s$ times the volume of
$W$ for some absolute positive constant $s$?

Here, by ``simple'' we mean a set which is a union of finitely many
dyadic cubes.

We prove that an affirmative answer to this question would have deep
consequences for the important space $BMO$ of functions of bounded
mean oscillation introduced by John and Nirenberg.

We recall and use the notion of a \textit{John-Str\"omberg pair }which
is closely related to the above question. The above mentioned result
is obtained as a consequence of a general result about these pairs.
We also present a number of additional results about these pairs.
\end{abstract}

\maketitle

\section{\label{sec:Intro}Introduction}

We are going to pose and discuss a rather simply formulated geometrical
and almost combinatorical question. We shall refer to it as \textbf{Question
A(1/2)}. Although it may seem to have no connection with functional
analysis, this question in fact has its origins in the study of a
very important space of functions. The space $BMO$, or $BMO(D)$,
which consists of all functions of bounded mean oscillation on a suitable
subset $D$ of $\mathbb{R}^{d}$, plays a surprisingly wide range
of different important roles in several branches of analysis. (For
some details about those roles see, e.g., \cite[p.~132]{CwikelMSagherYShvartsmanP2012}
or \cite[p.~4]{CwikelMSagherYShvartsmanP2010}.)

One of our aims here is to warmly invite readers, including those
with interests quite unrelated to spaces of functions or even to analysis
in general, to study our Question A(1/2), and hopefully even resolve
it. In particular (as our title seeks to indicate, and as confirmed
by \cite{HolzmanR2016}) those readers with expertise in geometry
in $\mathbb{R}^{n}$ or combinatorical problems in $\mathbb{Z}^{n}$
may well have some valuable insights. Those who choose to respond
to this invitation will be able to effectively do so, without needing
any familiarity whatsoever with the space $BMO$, nor with the proofs
in this paper, nor with the contents of the earlier papers which gave
rise to it.

Before formulating Question A(1/2) we shall mention another question
which motivates it and then fix some necessary (but rather simple)
terminology and notation.

This paper is a sequel to \cite{CwikelMSagherYShvartsmanP2012} and
its preliminary more detailed version \cite{CwikelMSagherYShvartsmanP2010}.
The main motivation for those papers, and also for this one, is the
wish for a better understanding of the celebrated John-Nirenberg inequality
\cite{JohnNirenberg}. This inequality is satisfied by every function
in the above mentioned space $BMO(D)$.

In particular, it is hoped to ultimately find an answer to:

\textbf{Question J-N.}\textit{ Can the constants in the John-Nirenberg
inequality for BMO functions of $d$ variables be chosen to be independent
of $d$?}

We refer to \cite{LernerA2013}, \cite{SlavinL2015}, \cite{SlavinLVasyninV2015},
\cite{SlavinLVasyuninV2016}, \cite{SlavinLZatitskiiP2019} and \cite{StolyarovDZatitskiyP2019},
and also to references in these papers. for some rather recent results
concerning best constants in certain versions of the John-Nirenberg
inequality. These deal mainly with the case where $d=1$. The results
of \cite{SlavinLZatitskiiP2019} apply to all values of $d$. Some
remarks on pp.~7--8 of \cite{SlavinL2015} recall some reasons for
being interested in the sizes of these constants, also for $d>1$.
We are grateful to Andrei Lerner, Pavel Shvartsman and Leonid Slavin
for information about these and related results.

As will be seen in a moment, our Question A(1/2) is quite easy to
formulate, and is expressed in terms readily accessible to a general
mathematical audience. Its importance lies in the fact that an affirmative
answer to it would imply an affirmative answer to Question J-N. It
would seem to be at least slightly easier to answer Question A(1/2)
than to answer certain very similar questions which have an analogous
role and which are posed in \cite{CwikelMSagherYShvartsmanP2010}
and \cite{CwikelMSagherYShvartsmanP2012}, and are also discussed
in the brief survey article \cite{CwikelMSagherYShvartsmanP2011}.
In the formulation of Question A(1/2), and indeed throughout this
paper, we shall use the following notation and terminology (much of
which is standard, and most of which was also used in \cite{CwikelMSagherYShvartsmanP2010}
and \cite{CwikelMSagherYShvartsmanP2012}).
\begin{defn}
\label{def:BasicStuff}We shall understand that

(i) $d$ is a positive integer, and that

(ii) the word \textbf{\textit{cube}} means a closed cube in $\mathbb{R}^{d}$
with sides parallel to the axes, i.e., the cartesian product of $d$
bounded closed intervals, all of the same length. Special roles will
be played by the $d$-dimensional \textbf{\textit{unit cube}} $\left[0,1\right]^{d}$
and its \textbf{\textit{dyadic subcubes}} (i.e, the cubes of the form
$\prod$$_{j=1}^{d}\left[(n_{j}-1)2^{-k},n_{j}2^{-k}\right]$, where
$k\in\mathbb{N}$ and the integers $n_{j}$ all satisfy $1\le n_{j}\le2^{k}$).
It will sometimes be convenient to use the standard notation $Q(x,r)$
for the cube of side length $2r$ centred at $x\in\mathbb{R}^{d}$,
i.e., $Q(x,r)=\left\{ y\in\mathbb{R}^{d}:\left\Vert y-x\right\Vert _{\ell_{d}^{\infty}}\le r\right\} $.

Furthermore,

(iii) we shall let $\mathcal{Q}(\mathbb{R}^{d})$ denote the collection
of all cubes in $\mathbb{R}^{d}$, and

(iv) the word $d$\textbf{\textit{-multi-cube}} will mean a non-empty
subset of $\left[0,1\right]^{d}$ which is the union of finitely many
\textit{dyadic} subcubes of $\left[0,1\right]^{d}$, and

(v) we will denote the $d$-dimensional Lebesgue measure of any Lebesgue
measurable subset $E$ of $\mathbb{R}^{d}$ by $\lambda(E)$. \textit{(But,
as explained below, for Question A(1/2) you do not need to know anything
at all about Lebesgue measure.)}

(The values of $d$ in (ii) and in (v) will be clear from the contexts
in which they appear.)
\end{defn}

In the formulation of Question A(1/2) we will not need any subtle
properties of Lebesgue measure. In fact we will only be dealing with
the simpler notion of $d$-dimensional volume. This is because the
only Lebesgue measurable subsets $E$ of $\mathbb{R}^{d}$ that we
will encounter in Question A(1/2) are cubes or $d$-multi-cubes and
their intersections with other cubes. All of these are finite unions
$E=\bigcup_{m=1}^{M}E_{m}$ of non-overlapping sets $E_{m}$, where
each $E_{m}$ is the cartesian product of $d$ closed intervals. So
we simply have that each $\lambda(E_{m})$ is the product of the lengths
of those $d$ intervals, and that $\lambda(E)=\sum_{m=1}^{M}\lambda(E_{m})$.

With some small adjustments, we can replace consideration of $\lambda(E)$
for all the above mentioned sets $E$ by consideration of the cardinality
of $E\cap2^{-k}\mathbb{Z}^{d}$ for a suitably large choice of $k\in\mathbb{N}$.
(This possibility was already very briefly hinted at in \cite{CwikelMSagherYShvartsmanP2011}.)
This is one reason for our suggestion above, in the opening paragraph
of this paper, that Question A(1/2) might perhaps also be approachable
via combinatorical considerations.

\smallskip{}

After these preparations, here at last is the promised question:

\smallskip{}

\textbf{Question A(1/2).}\textit{ Does there exist an absolute constant
$s>0$ which has the following property?}

\textbf{\textit{For every positive integer $d$, whenever $E_{+}$
and $E_{-}$ are two disjoint $d$-multi-cubes which satisfy 
\begin{equation}
\min\left\{ \lambda(E_{+}),\lambda(E_{-})\right\} >\frac{1}{2}\left(1-\lambda(E_{+})-\lambda(E_{-})\right)\,,\label{eq:Grink}
\end{equation}
 then there exists a cube $W$ which is contained in $[0,1]^{d}$
and for which
\begin{equation}
\min\left\{ \lambda(W\cap E_{+}),\lambda(W\cap E_{-})\right\} \ge s\lambda(W)\,.\label{eq:blurp}
\end{equation}
}}

Most of the work which is required to show that an affirmative answer
to Question A(1/2) implies an affirmative answer to Question J-N has
already been done in \cite{CwikelMSagherYShvartsmanP2010} and \cite{CwikelMSagherYShvartsmanP2012}.
Because of that we will not have any need at all here to deal with
any details concerning any of the versions of the function space $BMO$,
nor to even recall their definitions. It was shown in \cite{CwikelMSagherYShvartsmanP2010}
and \cite{CwikelMSagherYShvartsmanP2012} that an affirmative answer
to Question J-N would be a consequence of an affirmative answer to
the following question which was formulated at the beginnings of both
of those papers, and is clearly quite similar to Question A(1/2).

\textbf{Question A.}\textit{ Do there exist two absolute constants
$\tau\in(0,1/2)$ and $s>0$ which have the following property?}

\textbf{\textit{For every positive integer $d$ and for every cube
$Q$ in $\mathbb{R}^{d}$, whenever $E_{+}$ and $E_{-}$ are two
disjoint measurable subsets of $Q$ whose $d$-dimensional Lebesgue
measures satisfy 
\[
\min\left\{ \lambda(E_{+}),\lambda(E_{-})\right\} >\tau\lambda(Q\setminus E_{+}\setminus E_{-})\,,
\]
 then there exists a cube $W$ which is contained in $Q$ and for
which
\[
\min\left\{ \lambda(W\cap E_{+}),\lambda(W\cap E_{-})\right\} \ge s\lambda(W)\,.
\]
}}

It is obvious that an affirmative answer to Question A would imply
an affirmative answer to Question A(1/2). Our task here is to show
that the reverse implication also holds, namely that affirmatively
answering the apparently at least slightly less demanding Question
A(1/2) would suffice to also affirmatively answer Question A and therefore
also Question J-N.

In the next section we shall discuss some known results related to
Question A(1/2) including some limitations which can be anticipated
in any attempts to solve it. Then, after recalling some more notions
and providing some preliminary results in Section \ref{sec:Prelim}.
we will obtain the above mentioned reverse implication in Section
\ref{sec:Main}, as an immediate consequence (see Corollary \ref{cor:A-Half}
below) of the main result of this paper, Theorem \ref{thm:SmallerTau},
which is formulated in a slightly more general context than Questions
A and A(1/2).

We note that Theorem \ref{thm:SmallerTau} could also conceivably
be used to deduce an affirmative answer for Question J-N from certain
variants of Question A, which might perhaps be easier to answer than
Question A(1/2). Some explorations of the possibility of such options,
including some relevant numerical experiments, will perhaps be discussed
in a future sequel to this paper. The reader may perhaps care to take
note of some issues raised in Section 10 of \cite{CwikelMSagherYShvartsmanP2010}
which might turn out to be relevant for resolving Question A(1/2)
or its variants.

The above mentioned more general context in which Theorem \ref{thm:SmallerTau}
is formulated revolves around the notion of John-Str\"omberg pairs,
whose definition we will recall in Section \ref{sec:Prelim}. In Section
\ref{sec:Further} we shall present a number of additional results
about these pairs, beyond those which will be needed for our main
result.

Finally, the extremely brief Section \ref{sec:CSS-Comments} will
make some minor comments about the papers \cite{CwikelMSagherYShvartsmanP2010}
and \cite{CwikelMSagherYShvartsmanP2012}.

\section{\label{sec:KnownResults}Some constraints and some known results
related to Question A(1/2)}

As soon as we plunge into any attempt to answer Question A(1/2) we
can almost immediately identify our real ``enemies''. We can of
course always assume that $\lambda(E_{-})\le\lambda(E_{+})$. If we
restrict our attention to the case where $r_{0}\lambda(E_{+})\le\lambda(E_{-})\le\lambda(E_{+})$
for some fixed $r_{0}\in(0,1]$ then an obvious calculation shows
that, by simply choosing $W=[0,1]^{d}$, we can obtain \eqref{eq:blurp}
for $s=r_{0}/\left(3r_{0}+1\right)$. This shows that our above mentioned
``enemies'' are the cases where the ratio $\lambda(E_{-})/\lambda(E_{+})$
is arbitrarily small.

It is known that if Question A(1/2) has an affirmative answer, then
the positive number $s$ for which that answer holds must satisfy
\begin{equation}
s\le\sqrt{5}-2.\label{eq:smsa}
\end{equation}
This follows from the third of three results which have been obtained
by Ron Holzman \cite{HolzmanR2016} in connection with Question A(1/2).
It is a pleasure to describe these results. The first of them is an
affirmative answer to a non-trivial special case of Question A(1/2):
Holzman has shown, for every $d\in\mathbb{N}$, that whenever $E_{+}$
and $E_{-}$ are $d$-multi-cubes which are each finite unions of
dyadic cubes, \textit{all of side length} $1/2$, and they satisfy
\eqref{eq:Grink}, then there exists a cube $W$ in $\left[0,1\right]^{d}$
which satisfies \eqref{eq:blurp} for $s=1/4$. In fact, as will be
explained below in Remark \ref{rem:RonIsOptimal}, our Lemma \ref{lem:NotInJSd}
shows that this value of $s$ cannot be improved in the following
sense: If the general form of Question A(1/2) has an affirmative answer,
then the positive number $s$ which appears in that answer must satisfy
$s\le1/4.$ So this already comes quite close to establishing \eqref{eq:smsa}.

The second of Holzman's recent results is that this largest possible
value $1/4$ for $s$ is indeed attained in another special case of
Question A(1/2), where $d$ is restricted to take only one value,
namely $d=1$, but where there is no restriction on the form of the
disjoint measurable subsets $E_{+}$ and $E_{-}$ of $\left[0,1\right]$.
(See also Remark \ref{rem:RonAgain}.) Holzman's third result (see
Remark \ref{rem:RonIsOptimal} for details) uses a particular example
when $d=2$ to show that the analogue of his second result, when the
single chosen value for $d$ is greater than $1$, does not hold.
This example shows that for each $d\ge2$, as already stated above,
the relevant value of $s$ cannot be greater than $\sqrt{5}-2$.

\section{\label{sec:Prelim}Some further definitions and preliminary results}

The following definition (which is effectively the same as Definition
7.10 of \cite[p.~29]{CwikelMSagherYShvartsmanP2010} and Definition
7.9 of \cite[pp.~153--154]{CwikelMSagherYShvartsmanP2012}) recalls
a notion which plays a central role in \cite{CwikelMSagherYShvartsmanP2010}
and \cite{CwikelMSagherYShvartsmanP2012}. and which is of course
closely related to Question A.
\begin{defn}
\label{def:JohnStrombergPair} Let $d$ be a positive integer and
let $\mathcal{E}$ be a non-empty collection of Lebesgue measurable
subsets of $\mathbb{R}^{d}$. Suppose that each $E\in\mathcal{E}$
satisfies $0<\lambda(E)<\infty$. Let $\tau$ and $s$ be positive
numbers with the following property:

\textbf{({*})}\textbf{\textit{ Let $Q$ be an arbitrary set in $\mathcal{E}$
and let $E_{+}$ and $E_{-}$ be arbitrary disjoint measurable subsets
of $Q$. Suppose that
\begin{equation}
\min\left\{ \lambda(E_{+}),\lambda(E_{-})\right\} >\tau\lambda(Q\setminus E_{+}\setminus E_{-})\,.\label{eq:isms}
\end{equation}
 Then there exists a set $W\subset Q$ which is also in $\mathcal{E}$
and for which
\begin{equation}
\min\left\{ \lambda(E_{+}\cap W),\lambda(E_{-}\cap W)\right\} \ge s\lambda(W)\,.\label{eq:blerk}
\end{equation}
}} Then we will say that \textit{$\left(\tau,s\right)$ is a }\textbf{\textit{John-Strömberg
pair}}\textit{ for $\mathcal{E}$}.
\end{defn}

\begin{rem}
\label{rem:s-half}In the context of the preceding definition, if
$E_{|+}$, $E_{-}$ and $W$ are measurable subsets of $Q$ which
satisfy $E_{+}\cap E_{-}=\emptyset$ and $\lambda(W)>0$ and \eqref{eq:blerk}
for some $s>0$, then 
\[
\lambda(W)\ge\lambda(E_{+}\cap W)+\lambda(E_{-}\cap W)\ge2\min\left\{ \lambda(E_{+}\cap W),\lambda(E_{-}\cap W)\right\} \ge2s\lambda(W)>0
\]
and therefore $s\le1/2$. Consequently, for any choice of $\mathcal{E},$
any pair $(\tau,s)$ which is a John-Str\"omberg pair for $\mathcal{E}$
must satisfy $0<s\le1/2$. (Cf.~Remark 7.16 of \cite[p.~30]{CwikelMSagherYShvartsmanP2010}.)
\end{rem}

Several ``natural'' choices of the collection $\mathcal{E}$ are
mentioned on pp.~4--5 of \cite{CwikelMSagherYShvartsmanP2010} and
pp.~132--133 of \cite{CwikelMSagherYShvartsmanP2012}. These are
relevant for studying a number of variants of the function space $BMO$
which have been considered in the literature. The most ``classical''
of these choices, motivated by the context of \cite{JohnNirenberg}
and by Question J-N, is when we take $\mathcal{E}$ to be the collection
$\mathcal{Q}(\mathbb{R}^{d})$ of all cubes in $\mathbb{R}^{d}$ (as
introduced in Definition \ref{def:BasicStuff}(iii)) for some fixed
$d\in\mathbb{N}.$

It is known that $\left(\sqrt{2}-1,2^{-d}(3-2\sqrt{2})\right)$ is
a John-Strömberg pair for $\mathcal{Q}(\mathbb{R}^{d})$. (See Example
7.12(i) on p.~30 of \cite{CwikelMSagherYShvartsmanP2010}.) Indeed
in this paper we will need to use Definition \ref{def:JohnStrombergPair}
only in the special case where $\mathcal{E}$ is the collection $\mathcal{Q}(\mathbb{R}^{d})$.
But we have indicated the possibility of considering other choices
of the collection $\mathcal{E}$, since there may turn out to be variants
of our main result in this paper, and these may turn out to be relevant
for future research in the contexts of at least some of those other
choices. (We refer to \cite{CwikelMSagherYShvartsmanP2010} and \cite{CwikelMSagherYShvartsmanP2012}
for various results, in particular Theorem 9.1 of \cite[p.~41]{CwikelMSagherYShvartsmanP2010}
and \cite[p.~164]{CwikelMSagherYShvartsmanP2012}, which deal with
John-Strömberg pairs and their interactions with versions of the John-Nirenberg
inequality for variants of the space $BMO$ corresponding to other
choices of $\mathcal{E}$.)
\begin{defn}
\label{def:JSd}It will sometimes be convenient, for each $d\in\mathbb{N}$,
to let $JS(d)$ denote the set of all ordered pairs $\left(\tau,s\right)$
of positive real numbers which are John-Str\"omberg pairs for $\mathcal{Q}(\mathbb{R}^{d})$.
\end{defn}

Obviously, an equivalent reformulation of Question A is:

\textbf{Question A($\tau,s$).}\textit{ Do there exist two absolute
constants $\tau\in(0,1/2)$ and $s>0$ for which $\left(\tau,s\right)\in JS(d)$
for every $d\in\mathbb{N}$?}

This makes it very relevant to use the following known result. (Note
that below in Subsection \ref{subsec:Another} we will present a result
containing some new slight variants of it.)
\begin{thm}
\label{thm:onlyneedcubes}Let $d$ be a positive integer. The ordered
pair of positive numbers $\left(\tau,s\right)$ is a John-Str\"omberg
pair for $\mathcal{Q}(\mathbb{R}^{d})$ if and only if it has the
following property:

$(*)$ Whenever $F_{+}$ and $F_{-}$ are disjoint $d$-multi-cubes
which satisfy \textbf{
\begin{equation}
\min\left\{ \lambda(F_{+}),\lambda(F_{-})\right\} >\tau\lambda(\left[0,1\right]^{d}\setminus F_{+}\setminus F_{-})\,,\label{eq:Strict}
\end{equation}
}then there exists a cube $W$ contained in $[0,1]^{d}$ for which\textbf{
\begin{equation}
\min\left\{ \lambda(W\cap F_{+}),\lambda(W\cap F_{-})\right\} \ge s\lambda(W)\,.\label{eq:nezt}
\end{equation}
}
\end{thm}

\noindent \textit{Proof.} Obviously every John-Str\"omberg pair
$\left(\tau,s\right)$ for $\mathcal{Q}(\mathbb{R}^{d})$ must have
the property $(*)$. The fact that property $(*)$ implies that $\left(\tau,s\right)$
is a John-Str\"omberg pair for $\mathcal{Q}(\mathbb{R}^{d})$ is
exactly the content of Theorem 10.2 of \cite[p.~43]{CwikelMSagherYShvartsmanP2010},
and is proved on pp.~44--47 of \cite{CwikelMSagherYShvartsmanP2010}.
However, some other small issues remain to be clarified: It should
be mentioned that there are a few minor misprints in the proof of
Theorem 10.2 of \cite{CwikelMSagherYShvartsmanP2010}. But they do
not effect its validity. On page 46, in the third and fourth lines
above the inequality (10.16)) the notation $F_{\varepsilon}$ should
of course be changed to $F_{\flat}$ (in two places) and $\Omega_{\varepsilon}$
should be changed to $\Omega_{\flat}$. Then on line 11 of page 47
$r_{n}^{d}$ should be $\left(2r_{n}\right)^{d}$ and $r_{*}^{d}$
should be $\left(2r_{*}\right)^{d}$ . Similarly, four factors of
$2^{d}$ have been (harmlessly) omitted on line 5 of page 63 of \cite{CwikelMSagherYShvartsmanP2010}
in the proof of Lemma 10.1, which is an ingredient in the proof of
Theorem 10.2. That line should be
\[
=\left(2R_{n}\right)^{d}-\left(2r_{*}\right)^{d}+\left(2R_{n}\right)^{d}-\left(2r_{n}\right)^{d}.
\]
It should also be mentioned that, in the formulation of Theorem 10.2
of \cite{CwikelMSagherYShvartsmanP2010}, it is stated that the number
$\tau$ must satisfy $0<\tau<1/2$. But the proof is valid for all
$\tau>0$. (The requirement that $\tau<1/2$ was imposed only because
this is relevant in Theorem 9.1 of \cite[p.~41]{CwikelMSagherYShvartsmanP2010}
and \cite[p.~164]{CwikelMSagherYShvartsmanP2012}.)

The clarification of these small issues completes the proof of Theorem
\ref{thm:onlyneedcubes}. $\qed$

\section{\label{sec:Main}The main result}

We can now present our main result, and its obvious corollary for
dealing with Question A(1/2),
\begin{thm}
\label{thm:SmallerTau}Let $d$ be a positive integer. Suppose that
$(\tau,s)$ is a John-Str\"omberg pair for $\mathcal{Q}(\mathbb{R}^{d})$.
Then 
\begin{equation}
0<s\le1/2\label{eq:sLessThanHalf}
\end{equation}
 and, for each $\theta\in(0,s/(1-s))$, the pair $\left((1-\theta)\tau,s-\theta(1-s)\right)$
is also a John-Str\"omberg pair for $\mathcal{Q}(\mathbb{R}^{d})$.
\end{thm}

\begin{cor}
\label{cor:A-Half}Suppose that the answer to Question A(1/2) is affirmative,
i.e., suppose that $\left(1/2,s\right)$ is a John-Str\"omberg pair
for $\mathcal{Q}(\mathbb{R}^{d})$ for some positive number $s$ which
does not depend on $d$. Then, if we choose $\theta=s/(2-2s)$ in
Theorem \ref{thm:SmallerTau}, we obtain that $\left(\frac{1}{2}\cdot\frac{2-3s}{2-2s},\frac{s}{2}\right)$
is a John-Str\"omberg pair for $\mathcal{Q}\left(\mathbb{R}^{d}\right)$
for all $d\in\mathbb{N}$, which gives us an affirmative answer to
Question A$(\tau,s)$ and therefore also to Question A and Question
J-N.
\end{cor}

\noindent\textit{Proof of Theorem \ref{thm:SmallerTau}.} Let us
explicitly choose particular (necessarily positive) numbers $\tau$,
$s$ and $\theta$ which satisfy the hypotheses of the theorem.

The inequalities \eqref{eq:sLessThanHalf}, as observed in Remark
\ref{rem:s-half}, follow from the fact that $(\tau,s)$ is a John-Str\"omberg
pair for $\mathcal{Q}(\mathbb{R}^{d})$. Note that \eqref{eq:sLessThanHalf}
and the conditions imposed on $\theta$ ensure that $0<\theta<1$
and also that $0<\theta(1-s)<s$. Consequently, the numbers $\left(1-\theta\right)\tau$
and $s-\theta(1-s)$ are both strictly positive (as indeed they must
be if they are to form a John-Str\"omberg pair).

Throughout this proof we will let $Q$ denote the $d$-dimensional
unit cube, $Q=\left[0,1\right]^{d}$. Let $F_{+}$ and $F_{-}$ be
two arbitrary disjoint subsets of $Q$ which are both $d$-multi-cubes
and satisfy 
\begin{equation}
\min\left\{ \lambda(F_{+}),\lambda(F_{-})\right\} >(1-\theta)\tau\lambda\left(Q\setminus F_{+}\setminus F_{-}\right)\label{eq:fza}
\end{equation}
In view of Theorem \ref{thm:onlyneedcubes}, it will suffice to show
that there exists a subcube $W$ of $[0,1]^{d}$ such that 
\begin{equation}
\min\left\{ \lambda(W\cap F_{+}),\lambda(W\cap F_{-})\right\} \ge\left(s-\theta(1-s)\right)\lambda(W)\,.\label{eq:wwnt}
\end{equation}

For some sufficiently large integer $N$, both of the sets $F_{+}$
and $F_{-}$ are finite unions of dyadic cubes all of the same side
length $2^{-N}$, as of course is the whole unit cube $Q$. Since
the distance between $F_{+}$ and $F_{-}$ must be positive, there
must be some dyadic subcubes of $Q$ of side length $2^{-N}$ which
are not contained in $F_{+}\cup F_{-}$, and thus their interiors
are all contained in $Q\setminus F_{+}\setminus F_{-}$. Therefore
the set $Q\setminus F_{+}\setminus F_{-}$ is also, at least to within
some set of measure zero, a finite union of dyadic cubes of side length
$2^{-N}$.

Let $\mathcal{F}_{+}$ denote the collection of all dyadic cubes of
side length $2^{-N}$ which are contained in $F_{+}$. For each $\delta\in(0,2^{-N-1})$
and for each cube $W=\prod_{j=1}^{n}[a_{j},a_{j}+2^{-N}]$ in the
collection $\mathcal{F}_{+}$, let $H(W,\delta)$ denote the cube
concentric with $W$ contained in the interior of $W$ whose side
length is $2^{-N}-2\delta$, i.e., $H(W,\delta)=\prod_{j=1}^{n}[a_{j}+\delta,a_{j}+2^{-N}-\delta]$.
We will later use the following obvious fact: If $V$ is a cube which
intersects with $H(W,\delta)$ and has side length less than $\delta$,
then 
\begin{equation}
V\subset W\,.\label{eq:hpd}
\end{equation}
We now introduce a special subset $H_{\delta}^{+}$ of $F_{+}$ which
is ``slightly smaller'' than $F_{+}$. It is defined by
\[
H_{\delta}^{+}:=\bigcup_{W\in\mathcal{F}_{+}}H(W,\delta).
\]
With the help of this set, we will now choose a particular suitable
value for $\delta$ which must remain unchanged for the rest of this
proof. Obviously $\lambda(F_{+})-\lambda(H_{\delta}^{+})=\lambda(F_{+}\setminus H_{\delta}^{+})$
for each $\delta\in\left(0,2^{-N-1}\right)$ and $\lim_{\delta\searrow0}\lambda(F_{+}\setminus H_{\delta}^{+})=0$.
Therefore, since the inequality in \eqref{eq:fza} is strict, we can
and will choose our fixed value of $\delta$ to be sufficiently small
to ensure that 
\begin{equation}
\min\left\{ \lambda(H_{\delta}^{+}),\lambda(F_{-})\right\} >(1-\theta)\tau\lambda\left(Q\setminus F_{+}\setminus F_{-}\right)+\tau\lambda\left(F_{+}\setminus H_{\delta}^{+}\right)\,.\label{eq:zpe}
\end{equation}
Let $\mathcal{G}$ denote the collection of all dyadic cubes of side
length $2^{-N}$ whose interiors are contained in $Q\setminus F_{+}\setminus F_{-}$.
Given an arbitrary positive integer $k$, we divide each cube in the
collection $\mathcal{G}$ into $2^{dk}$ dyadic subcubes of side length
$2^{-N-k}$. Let $\mathcal{G}_{k}$ denote the collection of all dyadic
cubes obtained in this way. In other words, $\mathcal{G}_{k}$ is
simply the collection of all dyadic cubes of side length $2^{-N-k}$
whose interiors are contained in $Q\setminus F_{+}\setminus F_{-}$.
For each $W\in\mathcal{G}_{k}$ we let $U(W,k)$ be the (closed) cube
concentric with $W$ whose volume satisfies 
\begin{equation}
\lambda\left(U(W,k)\right)=\theta\lambda(W).\label{eq:Wtheta}
\end{equation}
We now introduce two more sets which will play particularly useful
roles for us. First we define
\[
V_{k}:=\bigcup_{W\in\mathcal{G}_{k}}U(W,k)
\]
and then use $V_{k}$ to define the disjoint union
\[
H_{k}^{-}:=F_{-}\cup V_{k}
\]
which contains and will be ``controllably'' larger than $F_{-}$.

In accordance with usual very standard notation, we will denote the
interior and the boundary of any given cube $W$ by $W^{\circ}$ and
$\partial W$ respectively. We of course have 
\begin{equation}
Q\setminus F_{+}\setminus F_{-}=\left(\bigcup_{W\in\mathcal{G}_{k}}W^{\circ}\right)\cup Z\label{eq:tvm}
\end{equation}
for some set $Z\subset\bigcup_{W\in\mathcal{G}_{k}}\partial W$. Since
$\lambda(Z)=0$, this implies that 
\begin{equation}
\lambda\left(Q\setminus F_{+}\setminus F_{-}\right)=\sum_{W\in\mathcal{G}_{k}}\lambda(W)\,.\label{eq:2ggp}
\end{equation}

In the following calculation we shall use \eqref{eq:tvm} in the second
line, and then, in the third line, the fact that, for every cube $W\in\mathcal{G}_{k}$,
the set $Z$ of zero measure introduced in \eqref{eq:tvm} is disjoint
from the cube $W^{\circ}$ and therefore also from $U(W,k)$. The
final line of the calculation will use the facts that, for each pair
of distinct cubes $W$ and $W'$ in $\mathcal{G}_{k}$, we have $U(W',k)\subset(W')^{\circ}$
and also $W^{\circ}\cap(W')^{\circ}=\emptyset$ and therefore $W^{\circ}\cap U(W',k)=\emptyset$.

Thus we obtain that

\begin{eqnarray*}
Q\setminus F_{+}\setminus H_{k}^{-} & = & Q\setminus F_{+}\setminus F_{-}\setminus\left(\bigcup_{W\in\mathcal{G}_{k}}U(W,k)\right)\\
 & = & \left(Z\cup\left(\bigcup_{W\in\mathcal{G}_{k}}W^{\circ}\right)\right)\setminus\left(\bigcup_{W\in\mathcal{G}_{k}}U(W,k)\right)\\
 & = & Z\cup\left(\left(\bigcup_{W\in\mathcal{G}_{k}}W^{\circ}\right)\setminus\left(\bigcup_{W\in\mathcal{G}_{k}}U(W,k)\right)\right)\\
 & = & Z\cup\left(\bigcup_{W\in\mathcal{G}_{k}}\left(W^{\circ}\setminus U(W,k)\right)\right)\,.
\end{eqnarray*}
Since $\lambda(Z)=0$, the preceding equalities, together with \eqref{eq:Wtheta}
and \eqref{eq:2ggp}, imply that
\begin{eqnarray}
\lambda(Q\setminus F_{+}\setminus H_{k}^{-}) & = & \sum_{W\in\mathcal{G}_{k}}\lambda\left(W^{\circ}\setminus U(W,k)\right)=(1-\theta)\sum_{W\in\mathcal{G}_{k}}\lambda(W)\nonumber \\
 & = & (1-\theta)\lambda(Q\setminus F_{+}\setminus F_{-})\,.\label{eq:NewGgp}
\end{eqnarray}

Any point in the set $Q\setminus H_{\delta}^{+}\setminus H_{k}^{-}$
which is not in $F_{+}$ must be in $Q\setminus F_{+}\setminus H_{k}^{-}$.
This shows that 
\begin{equation}
Q\setminus H_{\delta}^{+}\setminus H_{k}^{-}\subset\left(Q\setminus F_{+}\setminus H_{k}^{-}\right)\cup\left(F_{+}\setminus H_{\delta}^{+}\right)\,.\label{eq:rzm}
\end{equation}

Since the two sets in parentheses on the right side of \eqref{eq:rzm}
are disjoint, it follows that 
\begin{equation}
\lambda\left(Q\setminus H_{\delta}^{+}\setminus H_{k}^{-}\right)\le\lambda\left(Q\setminus F_{+}\setminus H_{k}^{-}\right)+\lambda\left(F_{+}\setminus H_{\delta}^{+}\right)\,.\label{eq:rcp}
\end{equation}
Now we can first use the fact that $F_{-}\subset H_{k}^{-}$ and then
invoke \eqref{eq:zpe} followed by \eqref{eq:NewGgp} and then \eqref{eq:rcp},
to obtain that 
\begin{align}
 & \phantom{=}\min\left\{ \lambda(H_{\delta}^{+}),\lambda(H_{k}^{-})\right\} \nonumber \\
 & \ge\min\left\{ \lambda(H_{\delta}^{+}),\lambda(F_{-})\right\} \nonumber \\
 & >(1-\theta)\tau\lambda\left(Q\setminus F_{+}\setminus F_{-}\right)+\tau\lambda\left(F_{+}\setminus H_{\delta}^{+}\right)\label{eq:NewVTT}\\
 & =\tau\lambda(Q\setminus F_{+}\setminus H_{k}^{-})+\tau\lambda\left(F_{+}\setminus H_{\delta}^{+}\right)\nonumber \\
 & \ge\tau\lambda(Q\setminus H_{\delta}^{+}\setminus H_{k}^{-})\,.\nonumber 
\end{align}
Since $V_{k}\subset Q\setminus F_{+}\setminus F_{-}$ and $F_{+}\cap F_{-}=\emptyset$,
we have 
\begin{equation}
H_{k}^{-}\subset F_{-}\cup\left(Q\setminus F_{+}\setminus F_{-}\right)=Q\setminus F_{+}\,.\label{eq:bmd}
\end{equation}
Furthermore $H_{\delta}^{+}\subset F_{+}$ and so $H_{\delta}^{+}$
and $H_{k}^{-}$ are disjoint measurable subsets of $Q$. According
to our hypotheses, $\left(\tau,s\right)$ is a John-Str\"omberg pair
for $\mathcal{Q}(\mathbb{R}^{d})$. Therefore, for each $k\in\mathbb{N}$,
since the quantity in the first line of \eqref{eq:NewVTT} is strictly
larger than the quantity in its last line, we deduce that there must
exist some subcube $W_{k}$ of $Q$ which satisfies
\begin{equation}
\min\left\{ \lambda(W_{k}\cap H_{\delta}^{+}),\lambda(W_{k}\cap H_{k}^{-})\right\} \ge s\lambda(W_{k})\,.\label{eq:fzz}
\end{equation}
We now claim that the side length of $W_{k}$, which we can conveniently
write as $\left(\lambda(W_{k})\right)^{1/d}$, must satisfy 
\begin{equation}
\left(\lambda(W_{k})\right)^{1/d}\ge\delta\,.\label{eq:dbv}
\end{equation}
To show this we first observe that, since the cube $W_{k}$ intersects
with $H_{\delta}^{+}$, it must intersect with the cube $H(W,\delta)$
for at least one cube $W$ in the collection $\mathcal{F}_{+}$. If
\eqref{eq:dbv} does not hold, i.e., if $W_{k}$ has side length less
than $\delta$, then (cf.~the discussion immediately preceding \eqref{eq:hpd})
$W_{k}$ must be completely contained in that particular cube $W$
and therefore also in $F_{+}$. Consequently $W_{k}$ cannot intersect
with the set $H_{k}^{-}$. (Here we have used \eqref{eq:bmd} once
more.) This contradicts \eqref{eq:fzz} and shows that \eqref{eq:dbv}
does hold.

Let $\widetilde{W}_{k}$ be a new cube containing $W_{k}$ and concentric
with $W_{k}$. More precisely, if $W_{k}=\prod_{j=1}^{d}\left[a_{j},b_{j}\right]$,
then we choose $\widetilde{W}_{k}$ to be $\prod_{j=1}^{d}\left[a_{j}-2^{-N-k},b_{j}+2^{-N-k}\right]$.
Clearly $\widetilde{W}_{k}$ contains all dyadic cubes of side length
$2^{-N-k}$ which intersect with $W_{k}$, and the side length $\left(\lambda(\widetilde{W}_{k})\right)^{1/d}$
of $\widetilde{W}_{k}$ satisfies 
\[
\left(\lambda(\widetilde{W}_{k})\right)^{1/d}=\left(\lambda(W_{k})\right)^{1/d}+2^{1-N-k}\,.
\]
This, together with \eqref{eq:dbv}, gives us that 
\begin{eqnarray*}
\left(\frac{\lambda\left(\widetilde{W}_{k}\right)}{\lambda\left(W_{k}\right)}\right)^{1/d} & = & \frac{\left(\lambda\left(W_{k}\right)\right)^{1/d}+2^{1-N-k}}{\left(\lambda\left(W_{k}\right)\right)^{1/d}}\\
 & = & 1+\frac{2^{1-N-k}}{\left(\lambda\left(W_{k}\right)\right)^{1/d}}\le1+\frac{2^{1-N-k}}{\delta}\,.
\end{eqnarray*}
It follows that 
\begin{eqnarray}
\lambda\left(\widetilde{W}_{k}\right)-\lambda(W_{k}) & = & \lambda(W_{k})\left(\frac{\lambda\left(\widetilde{W}_{k}\right)}{\lambda\left(W_{k}\right)}-1\right)\nonumber \\
 & \le & \lambda(W_{k})\left(\left(1+\frac{2^{1-N-k}}{\delta}\right)^{d}-1\right)\,.\label{eq:osfd}
\end{eqnarray}
Let $\mathcal{H}_{k}$ be the collection of all cubes in $\mathcal{G}_{k}$
which intersect with $W_{k}$. Clearly 
\begin{equation}
\bigcup_{W\in\mathcal{H}_{k}}W^{\circ}\subset(Q\setminus F_{+}\setminus F_{-})\cap\widetilde{W}_{k}.\label{eq:GKT}
\end{equation}
 Our definitions of $V_{k}$ and of $\mathcal{H}_{k}$ also immediately
give us that 
\begin{eqnarray*}
W_{k}\cap V_{k} & = & W_{k}\cap\left(\bigcup_{W\in\mathcal{G}_{k}}U(W,k)\right)\\
 & = & W_{k}\cap\left(\bigcup_{W\in\mathcal{H}_{k}}U(W,k)\right)\\
 & \subset & \bigcup_{W\in\mathcal{H}_{k}}U(W,k)\,.
\end{eqnarray*}
This and then \eqref{eq:Wtheta} will enable us to obtain the first
line in the following calculation. Its second line will use the fact
that the interiors of the cubes in $\mathcal{G}_{k}$ and therefore
also in $\mathcal{H}_{k}$, are pairwise disjoint. Its third line
will use \eqref{eq:GKT}. Its fourth line will use the fact that $H_{\delta}^{+}\subset F_{+}$.
Its sixth line will follow from the obvious inclusion $W_{k}\subset\widetilde{W}_{k}$.
The justifications of all other steps should be evident. (We wonder,
casually, whether it might somehow be possible to replace the simple-minded
transition from the sixth to the seventh line by a sharper estimate,
which might then lead to a (slightly) stronger version of Theorem
\ref{thm:SmallerTau}.) 
\begin{eqnarray*}
\lambda(W_{k}\cap V_{k}) & \le & \sum_{W\in\mathcal{H}_{k}}\lambda(U(W,k))=\theta\sum_{W\in\mathcal{H}_{k}}\lambda(W)\\
 & = & \theta\sum_{W\in\mathcal{H}_{k}}\lambda(W^{\circ})=\theta\lambda\left(\bigcup_{W\in\mathcal{H}_{k}}W^{\circ}\right)\\
 & \le & \theta\lambda\left((Q\setminus F_{+}\setminus F_{-})\cap\widetilde{W}_{k}\right)\\
 & \le & \theta\lambda\left((Q\setminus H_{\delta}^{+}\setminus F_{-})\cap\widetilde{W}_{k}\right)\\
 & \le & \theta\lambda\left((Q\setminus H_{\delta}^{+}\setminus F_{-})\cap W_{k}\right)+\theta\lambda\left(\widetilde{W}_{k}\setminus W_{k}\right)\\
 & = & \theta\lambda\left((Q\setminus H_{\delta}^{+}\setminus F_{-})\cap W_{k}\right)+\theta\left(\lambda\left(\widetilde{W}_{k}\right)-\lambda\left(W_{k}\right)\right)\\
 & \le & \theta\lambda\left((Q\setminus H_{\delta}^{+})\cap W_{k}\right)+\theta\left(\lambda\left(\widetilde{W}_{k}\right)-\lambda\left(W_{k}\right)\right)\\
 & = & \theta\lambda\left(W_{k}\setminus H_{\delta}^{+}\right)+\theta\left(\lambda\left(\widetilde{W}_{k}\right)-\lambda\left(W_{k}\right)\right)\,.
\end{eqnarray*}

Combining the result of this calculation with \eqref{eq:osfd}, we
deduce that, for each $k\in\mathbb{N}$, 
\begin{equation}
\lambda(W_{k}\cap V_{k})\le\theta\lambda\left(W_{k}\setminus H_{\delta}^{+}\right)+\varepsilon_{k}\lambda(W_{k})\label{eq:ziuy}
\end{equation}
where $\varepsilon_{k}:=\theta\left(\left(1+\frac{2^{1-N-k}}{\delta}\right)^{d}-1\right)$
and we will later use the obvious fact that
\begin{equation}
\lim_{k\to\infty}\varepsilon_{k}=0\,.\label{eq:db2}
\end{equation}

In view of \eqref{eq:fzz}, we have that 
\begin{equation}
\lambda(W_{k}\setminus H_{\delta}^{+})=\lambda(W_{k})-\lambda(W_{k}\cap H_{\delta}^{+})\le(1-s)\lambda(W_{k})\,.\label{eq:dyz}
\end{equation}
We now have the ingredients needed to estimate $\lambda\left(W_{k}\cap F_{-}\right)$
from below. Since $\lambda(W_{k}\cap H_{k}^{-})=\lambda(W_{k}\cap F_{-})+\lambda(W_{k}\cap V_{k})$
we can use \eqref{eq:fzz} and \eqref{eq:ziuy} and then \eqref{eq:dyz}
to obtain that 
\begin{eqnarray}
\lambda(W_{k}\cap F_{-}) & = & \lambda(W_{k}\cap H_{k}^{-})-\lambda(W_{k}\cap V_{k})\nonumber \\
 & \ge & s\lambda(W_{k})-\theta\lambda\left(W_{k}\setminus H_{\delta}^{+}\right)-\varepsilon_{k}\lambda(W_{k})\nonumber \\
 & \ge & s\lambda(W_{k})-\theta\left(1-s\right)\lambda(W_{k})-\varepsilon_{k}\lambda(W_{k})\nonumber \\
 & = & (s-\theta(1-s)-\varepsilon_{k})\lambda(W_{k})\,.\label{eq:ltz}
\end{eqnarray}

Furthermore, again with the help of \eqref{eq:fzz}, we also obviously
have that 
\begin{equation}
\lambda(W_{k}\cap F_{+})\ge\lambda(W_{k}\cap H_{\delta}^{+})\ge s\lambda(W_{k})\,.\label{eq:xir}
\end{equation}

\smallskip{}

As already observed at the beginning of this proof, $s-\theta(1-s)$
and $\theta(1-s)$ are both strictly positive. So, in view of \eqref{eq:db2},
for all sufficiently large $k$, we will also have $s>(s-\theta(1-s)-\varepsilon_{k})$
and $(s-\theta(1-s)-\varepsilon_{k})>0$. Therefore, we can now deduce
from \eqref{eq:fza} and \eqref{eq:ltz} and \eqref{eq:xir} that
$\left((1-\theta)\tau,s-\theta(1-s)-\varepsilon_{k}\right)$ is a
John-Str\"omberg pair for $\mathcal{Q}(\mathbb{R}^{d})$ for all
sufficiently large $k$. This seems to be very close to our required
result. Indeed it will only require a little more effort to obtain
that result.

\smallskip{}

For each $k\in\mathbb{N}$, let $x_{k}$ be the centre of the cube
$W_{k}$ and let $r_{k}$ be half its side length. I.e., we can set
$W_{k}=Q(x_{k},r_{k})$ in the standard notation recalled above in
Definition \ref{def:BasicStuff}(ii). Of course $x_{k}\in Q$ and,
by \eqref{eq:dbv} and the fact that $W_{k}\subset Q$, we also have
$\delta/2\le r_{k}\le1/2$. Therefore there exists a strictly increasing
sequence $\left\{ n_{k}\right\} _{k\in\mathbb{N}}$ of positive integers
such that the sequences $\left\{ x_{n_{k}}\right\} _{k\in\mathbb{N}}$
and $\left\{ r_{n_{k}}\right\} _{k\in\mathbb{N}}$ converge, respectively,
to a point $x\in Q$ and a number $r\in[\delta/2,1/2]$. Let $W$
be the cube $W=Q(x,r)$. Then, by Lemma 10.1 of \cite[p.~43]{CwikelMSagherYShvartsmanP2010}
and \eqref{eq:ltz} and \eqref{eq:db2}, it follows that $W\subset Q$
and 
\begin{eqnarray*}
\lambda(W\cap F_{-}) & = & \lim_{k\to\infty}\lambda(W_{n_{k}}\cap F_{-})\\
 & \ge & \lim_{k\to\infty}(s-\theta(1-s)-\varepsilon_{n_{k}})\lambda(W_{n_{k}})\\
 & = & (s-\theta(1-s))\lambda(W)\,.
\end{eqnarray*}
Similarly, by replacing $k$ by $n_{k}$ in \eqref{eq:xir} and then
passing to the limit as $k$ tends to $\infty$, we obtain that 
\[
\lambda(W\cap F_{+})\ge s\lambda(W)\ge(s-\theta(1-s))\lambda(W)\,.
\]
All this shows that the subcube $W$ of $Q$ satisfies \eqref{eq:wwnt}
and therefore completes the proof of the theorem. $\qed$

\bigskip{}

\begin{rem}
\label{rem:MaybeImprove}It seems quite possible that a more elaborate
version of the preceding proof might show that the hypotheses of Theorem
\ref{thm:SmallerTau} can yield a stronger result, namely that $\left(\tau',s'\right)$
is a John-Str\"omberg pair for some positive number $\tau'$ smaller
than $\left(1-\theta\right)\tau$ and/or for some number $s'$ greater
than $s-\theta(1-s)$. Perhaps a strategy for proving this might involve
separately considering the two cases where $\lambda(F_{+})>c\lambda(F_{-})$
and where $\lambda(F_{-})\le\lambda(F_{+})\le c\lambda(F_{-})$ for
some suitably chosen constant $c>1$.
\end{rem}

\section{\label{sec:Further}Some further results about John-Str\"omberg
pairs}

\subsection{\label{subsec:Another}Other characterizations of John-Str\"omberg
pairs.}

It will be convenient to begin by introducing some ``technical''
terminology.
\begin{defn}
\label{def:Tame} Let $F_{+}$ and $F_{-}$ be two disjoint measurable
subsets of $\left[0,1\right]^{d}$ which both have positive measure
and which also satisfy one of the following two conditions:

\noindent (i) $\lambda\left([0,1]^{d}\setminus F_{+}\setminus F_{-}\right)=0$.

\noindent (ii) There exists a non-empty open subset $\Omega$ of $\left[0,1\right]^{d}$
for which

\noindent 
\begin{equation}
\max\left\{ \lambda(F_{+}\cap\Omega),\lambda(F_{-}\cap\Omega)\right\} <\lambda(\Omega)\,\text{\,\,\ and\,\,\,}\min\left\{ \lambda(F_{+}\cap\Omega),\lambda(F_{-}\cap\Omega)\right\} =0.\label{eq:WithOmega}
\end{equation}

Then we shall say that $\left(F_{+},F_{-}\right)$ is a \textbf{\textit{tame
couple}} in $\left[0,1\right]^{d}$.
\end{defn}

\begin{rem}
\label{rem:TwoEquivalentConditions} Obviously the two conditions
(i) and (ii) are mutually exclusive. Furthermore, the following simple
argument shows that the condition (ii) is equivalent to

(ii)' There exists a dyadic cube $V$ contained in $\left[0,1\right]^{d}$
for which

\noindent 
\begin{equation}
\max\left\{ \lambda(F_{+}\cap V),\lambda(F_{-}\cap V)\right\} <\lambda(V)\,\text{\,\,\ and\,\,\,}\min\left\{ \lambda(F_{+}\cap V),\lambda(F_{-}\cap V)\right\} =0.\label{eq:VinsteadOfOmega}
\end{equation}

The Lebesgue measures of each of the sets appearing in \eqref{eq:VinsteadOfOmega}
remain unchanged if the dyadic cube $V$ in them is replaced by its
interior $V^{\circ}$. Therefore condition (ii)' implies condition
(ii). Conversely, suppose that some open set $\Omega$ satisfies \eqref{eq:WithOmega}.We
can suppose that $0=\lambda(F_{-}\cap\Omega)\le\lambda(F_{+}\cap\Omega)<\lambda(\Omega)$,
(Otherwise simply reverse the roles of $F_{+}$ and $F_{-}$ in what
is to follow.) By Theorem 1.11 of \cite[p.~8]{WheedenRZygmundA1977},
$\Omega$ can be expressed as the union of a sequence of non-overlapping
dyadic cubes $\Omega=\bigcup_{n\in\mathbb{N}}V_{n}$. So $\sum_{n=1}^{\infty}\lambda(F_{+}\cap V_{n})<\sum_{n=1}^{\infty}\lambda(V_{n})$
and consequently $\lambda(F_{+}\cap V_{n})<\lambda(V_{n})$ for at
least one $n$ which of course also satisfies $\lambda(F_{-}\cap V_{n})=0$.
So, for that choice of $n$, the dyadic cube $V_{n}$ satisfies \eqref{eq:VinsteadOfOmega}.
Therefore condition (ii) implies condition (ii)' and so (ii) and (ii)'
are indeed equivalent.

\smallskip{}
\end{rem}

We can now present the new variant of Theorem \ref{thm:onlyneedcubes}
to which we referred in the preamble to that theorem. It will specify
two more properties of a pair $\left(\tau,s\right)$, which we shall
label as $(**)$ and $(***)$, and which are each equivalent to the
property that $\left(\tau,s\right)\in JS(d)$. Note that the only
difference between the statement of property $(**)$ in Theorem \ref{thm:EqualityAllowed},
and the statement of property $(*)$ in Theorem \ref{thm:onlyneedcubes}
is that the strict equality ``$>$'' which appears in \eqref{eq:Strict}
in the statement of $(*)$ has been replaced by ``$\ge$'' in \eqref{eq:EqualityAllowed}
in the statement of $(**)$. Property $(***)$ of Theorem \ref{thm:EqualityAllowed}
is more elaborate, and requires the terminology of Definition \ref{def:Tame}.
\begin{thm}
\label{thm:EqualityAllowed}Let $d$ be a positive integer and let
$\left(\tau,s\right)$ be an ordered pair of positive numbers. Then
each of the following two properties is equivalent to the property
that $\left(\tau,s\right)$ is a John-Str\"omberg pair for $\mathcal{Q}(\mathbb{R}^{d})$.

$(**)$ Whenever $F_{+}$ and $F_{-}$ are disjoint $d$-multi-cubes
which satisfy \textbf{
\begin{equation}
\min\left\{ \lambda(F_{+}),\lambda(F_{-})\right\} \ge\tau\lambda(\left[0,1\right]^{d}\setminus F_{+}\setminus F_{-})\,,\label{eq:EqualityAllowed}
\end{equation}
}then there exists a cube $W$ contained in $[0,1]^{d}$ for which\textbf{
\begin{equation}
\min\left\{ \lambda(W\cap F_{+}),\lambda(W\cap F_{-})\right\} \ge s\lambda(W)\,.\label{eq:Bloomp}
\end{equation}
$(***)$ }Whenever $\left(F_{+},F_{-}\right)$ is a tame couple of
measurable subsets of $\left[0,1\right]^{d}$ which satisfies \eqref{eq:EqualityAllowed},
then there exists a cube $W$ contained in $\left[0,1\right]^{d}$
for which \eqref{eq:Bloomp} holds.
\end{thm}

\noindent \textit{Proof.} Let us first show that property $(***)$
implies property $(**)$. Suppose that $F_{+}$ and $F_{-}$ are arbitrary
disjoint $d$-multi-cubes. By the same simple reasoning as was given
in the paragraph immediately following \eqref{eq:wwnt} in the proof
of Theorem \ref{thm:SmallerTau}, we see that the set $[0,1]^{d}\setminus F_{+}\setminus F_{-}$
must contain the interior $V^{\circ}$ of at least one dyadic subcube
$V$ of $\left[0,1\right]^{d}$. The set $\Omega=V^{\circ}$ of course
satisfies the condition (ii) in Definition \ref{def:Tame}. Consequently,
$\left(F_{+},F_{-}\right)$ is a tame couple in $\left[0.1\right]^{d}$.
It immediately follows that property $(***)$ indeed does imply property
$(**)$.

Obviously, if the pair $\left(\tau,s\right)$ of positive numbers
has property $(**)$, then it also has property $(*)$ of Theorem
\ref{thm:onlyneedcubes} and therefore $\left(\tau,s\right)$ is a
John-Str\"omberg pair for $\mathcal{Q}(\mathbb{R}^{d})$.

Now suppose that $\left(\tau,s\right)$ is a John-Str\"omberg pair
for $\mathcal{Q}(\mathbb{R}^{d})$. Then, by letting $\left\{ \tau_{n}\right\} _{n\in\mathbb{N}}$
and $\left\{ s_{n}\right\} _{n\in\mathbb{N}}$ be the special constant
sequences $\tau_{n}=\tau$ and $s_{n}=s$, we see that $\left(\tau,s\right)$
satisfies the hypotheses of Lemma \ref{lem:Sequences of JS pairs}
which we will formulate and prove below. In view of part (b) of that
lemma, $\left(\tau,s\right)$ has property $\left(***\right)$ . This
completes the proof of the theorem. $\qed$
\begin{rem}
It is very natural to wonder whether property $(***)$ is equivalent
to a stronger and more simply expressed variant of that property in
which the same implication is required to hold for \textit{all} pairs
of disjoint measurable subsets $F_{+}$, $F_{-}$ of $\left[0,1\right]^{d}$
which both have positive measure, thus omitting the requirement that
$F_{+}$ and $F_{-}$ should form a tame couple. (In other words,
as in Remark 7.17 of \cite[p.~30]{CwikelMSagherYShvartsmanP2010},
we are essentially wondering whether in Definition \ref{def:JohnStrombergPair},
at least in the case where $\mathcal{E}=\mathcal{Q}(\mathbb{R}^{d})$,
it would be equivalent to replace ``$>$'' by ``$\ge$'' in \eqref{eq:isms},
of course then with the necessary proviso that $\lambda(E_{+})$ and
$\lambda(E_{-})$ are both positive.) We are unable to answer this
question, but Lemma \ref{lem:Sequences of JS pairs} shows that if
its answer is negative, then the sets $F_{+}$ and $F_{-}$which provide
a counterexample, must both have very intricate structure (and it
would not be inappropriate to refer to them as forming a ``wild''
couple).
\end{rem}

\begin{lem}
\label{lem:Sequences of JS pairs} Let $d$ be a positive integer.
Let $\left\{ \tau_{n}\right\} _{n\in\mathbb{N}}$ and $\left\{ s_{n}\right\} _{n\in\mathbb{N}}$
be two sequences of positive numbers which converge to positive limits,
$\tau$ and $s$ respectively. Suppose that $\left(\tau_{n},s_{n}\right)$
is a John-Str\"omberg pair for $\mathcal{Q}(\mathbb{R}^{d})$ for
every $n\in\mathbb{N}$. Then

$(\mathrm{a})$ the limiting pair $\left(\tau,s\right)$ has the property
$(**)$ of Theorem \ref{thm:EqualityAllowed}, and,

$(\mathrm{b})$ in the special case where $s_{n}\ge s$ for every
$n\in\mathbb{N}$, the limiting pair $\left(\tau,s\right)$ also has
the property $(***)$ of Theorem \ref{thm:EqualityAllowed}.
\end{lem}

\noindent \textit{Proof.} In this proof it will be convenient to
use a trivial generalization of the standard notation for cubes recalled
in Definition \ref{def:BasicStuff}(ii) and to introduce some additional
obvious terminology as follows:
\begin{defn}
\label{def:LimitOfCubes} We permit ourselves to extend the notation
$Q(x,r)$ for a cube of side length $r$ centred at $x$, for each
$x\in\mathbb{R}^{d}$ and $r>0$, also to the case where $r=0$, by
letting $Q(x,0)$ denote the singleton $\left\{ x\right\} .$ Let
$\left\{ x_{k}\right\} _{k\in\mathbb{N}}$ be a convergent sequence
in $\mathbb{R}^{d}$ and let $\left\{ r_{k}\right\} _{k\in\mathbb{N}}$
be a convergent sequence of positive numbers. For each $k\in\mathbb{N}$
let $W_{k}$ be the cube $Q(x_{k},r_{k})$. Then we refer to $\left\{ W_{k}\right\} _{k\in\mathbb{N}}$
as a \textit{convergent sequence of cubes}, and define its limit $\lim_{k\to\infty}W_{k}$
to be the cube or singleton $Q\left(\lim_{k\to\infty}x_{k},\lim_{k\to\infty}r_{k}\right)$.
\end{defn}

Let us fix an arbitrary $d\in\mathbb{N}$, and arbitrary sequences
$\left\{ \tau_{n}\right\} _{n\in\mathbb{N}}$ and $\left\{ s_{n}\right\} _{n\in\mathbb{N}}$
of positive numbers which tend respectively to the positive numbers
$\tau$ and $s$, and have the property that $\left(\tau_{n},s_{n}\right)\in JS(d)$
for each $n\in\mathbb{N}$.

In order to deduce that $\left(\tau,s\right)$ has property $(**)$
and/or property $(***)$, we begin by fixing two arbitrary disjoint
measurable sets $F_{+}$ and $F_{-}$ which are contained in of $\left[0,1\right]^{d}$,
which form a tame couple in $\left[0,1\right]^{d}$, and which also
satisfy \eqref{eq:EqualityAllowed}.

To obtain part (a) of the lemma we have the task of proving, for these
choices of $F_{+},$ $F_{-}$, $\left\{ \tau_{n}\right\} _{n\in\mathbb{N}}$,
$\left\{ s_{n}\right\} _{n\in\mathbb{N}}$, $\tau$ and $s$, that
there exists a subcube $W$ of $\left[0,1\right]^{d}$ which satisfies
\eqref{eq:Bloomp} . If needed at any stage of our proof of that fact,
we may make the additional assumption that $F_{+}$ and $F_{-}$ are
both $d$-multi-cubes.

An analogous task is required to obtain part (b) of the lemma, with
the difference that in our proof this time, instead of being able
to assume that $F_{+}$ and $F_{-}$ are $d$-multi-cubes, we may
make the additional assumption, if needed, that $s_{n}\ge s$ for
all $n$.

There is quite a lot of overlap in the ingredients which will be used
for performing these two tasks, and it may help avoid some confusion
if we give some general description, in advance, of each of the four
steps which we shall use to accomplish both of them almost simultaneously.
We stress that neither of the two above mentioned additional assumptions
will be needed in the first two of these four steps.

In Step 1 of the proof, we shall use the given measurable sets $F_{+}$
and $F_{-}$ and the given sequences $\left\{ \tau_{n}\right\} _{n\in\mathbb{N}}$
and $\left\{ s_{n}\right\} _{n\in\mathbb{N}}$ to construct a special
sequence $\left\{ W(n_{k})\right\} _{k\in\mathbb{N}}$ of subcubes
of $\left[0,1\right]^{d}$ which converges either to a cube or to
a singleton set, in the sense of Definition \ref{def:LimitOfCubes}.
We shall perform this construction in two different (and unexplained
and sometimes complicated) ways, depending on which of the two conditions
of Definition \ref{def:Tame} is satisfied by $\left(F_{+},F_{-}\right)$.
It is only in later steps of the proof that we will be able to properly
see the usefulness of the particular features of these constructions.

In Step 2, we shall see that if the sequence $\left\{ W(n_{k})\right\} _{k\in\mathbb{N}}$
obtained in the preceding step converges to a cube $W$, then that
cube is contained in $\left[0,1\right]^{d}$ and satisfies \eqref{eq:Bloomp}.

In Step 3, we shall see that whenever the disjoint measurable sets
$F_{+}$ and $F_{-}$ are both required to also be $d$-multi-cubes,
then the sequence $\left\{ W(n_{k})\right\} _{k\in\mathbb{N}}$ can
be always be constructed so that it converges to a subcube rather
than a singleton. In view of Step 2, this will complete the proof
of part (a) of the lemma.

Finally, in Step 4, in order to complete the proof of part (b), it
will remain (again in view of Step 2) only to deal with the case where
the sequence $\left\{ W(n_{k})\right\} _{k\in\mathbb{N}}$ converges
to a singleton. We will do this by showing that, in this case, when
we also impose the requirement that the sequence $\left\{ s_{n}\right\} _{n\in\mathbb{N}}$
satisfies $s_{n}\ge s$ for each $n$, then there exists a integer
$k_{0}$ for which the subcube $W=W(n_{k_{0}})$ is contained in $\left[0,1\right]^{d}$
and satisfies \eqref{eq:Bloomp}.

\textbf{\textit{STEP 1: Construction of the special convergent sequence
$\left\{ W(n_{k})\right\} _{k\in\mathbb{N}}$.}}

Since $\left(F_{+},F_{-}\right)$ is a tame couple, it must satisfy
either condition (i) or condition (ii) of Definition \ref{def:Tame}.

Our construction of $\left\{ W(n_{k})\right\} _{k\in\mathbb{N}}$
will be quite simple in the case where condition (i) holds, i.e.,
when 
\begin{equation}
\lambda\left([0,1]^{d}\setminus F_{+}\setminus F_{-}\right)=0.\label{eq:NuStum}
\end{equation}
The strict inequalities $0<\lambda(F_{+})<\lambda\left(\left[0,1\right]^{d}\right)$
enable us to use Lemma 7.1 of \cite[p.~25]{CwikelMSagherYShvartsmanP2010}
or of \cite[p.~150]{CwikelMSagherYShvartsmanP2012} (whose simple
proof is a special case of the proof of Lemma 7.5 on pages 26-27 of
\cite{CwikelMSagherYShvartsmanP2010}) to ensure the existence of
a cube $W=Q(x,r)\subset[0,1]^{d}$ (of course with $r>0$) for which
\begin{equation}
\lambda(W\cap F_{+})=\lambda(W\setminus F_{+})=\frac{1}{2}\lambda(W).\label{eq:EasyStuff}
\end{equation}
The condition \eqref{eq:NuStum} is of course equivalent to $\lambda(F_{+}\cup F_{-})=1$
and to the fact that the sets $F_{-}$ and $\left[0,1\right]^{d}\setminus F_{+}$
coincide to with sets of measure zero. Therefore we also have
\begin{equation}
\lambda\left(W\cap F_{-}\right)=\lambda\left(W\setminus F_{+}\right).\label{eq:same}
\end{equation}
In this case we will contruct our required special convergent sequence
$\left\{ W(n_{k})\right\} _{k\in\mathbb{N}}$ by simply setting $x_{k}=x$,
$r_{k}=r$, $W(k)=W$ and $n_{k}=k$ for each $k\in\mathbb{N}$. Then
the sequence $\left\{ W(n_{k})\right\} _{k\in\mathbb{N}}$ of course
converges, not to a singleton, but to the cube $W$ which is its constant
value.

It remains to consider the more complicated case when $F_{+}$ and
$F_{-}$ satisfy condition (ii) or equivalently condition (ii)' of
Remark \ref{rem:TwoEquivalentConditions}.

In that case, since the roles of $F_{+}$ and $F_{-}$ are interchangeable
in properties $(**)$ and $(***)$ and in Definition \ref{def:Tame}
and in \eqref{eq:Bloomp}, we can assume, without loss of generality,
that, by Remark \ref{rem:TwoEquivalentConditions}, there exists a
dyadic cube $V$ contained in $\left[0,1\right]^{d}$ for which 
\begin{equation}
0=\lambda(F_{-}\cap V)\le\lambda(F_{+}\cap V)<\lambda(V).\label{eq:OmegaStuff}
\end{equation}
Since the interior $V^{\circ}$ of $V$ satisfies
\[
\lambda(V^{\circ}\setminus F_{+})=\lambda(V\setminus F_{+})=\lambda(V)-\lambda(V\cap F_{+})>0,
\]
the Lebesgue differentiation theorem guarantees the existence of a
point $z\in V^{\circ}\setminus F_{+}$ for which $\lim_{r\searrow0}\frac{\lambda\left(Q(z,r)\cap(V^{\circ}\setminus F_{+})\right)}{\lambda(Q(z,r))}=1$.
Since $Q(z,r)\subset V^{\circ}$ for all sufficiently small $r$,
we also have $\lim_{r\searrow0}\frac{\lambda\left(Q(z,r)\setminus F_{+})\right)}{\lambda(Q(z,r))}=1$
and, consequently, $\lim_{r\searrow0}\frac{\lambda\left(Q(z,r)\cap F_{+}\right)}{\lambda(Q(z,r))}=0$.
These properties of $V^{\circ}$ and $z$ enable us to assert the
existence of a sequence $\left\{ \rho_{k}\right\} _{k\in\mathbb{N}}$
of positive numbers tending to monotonically to zero, such that for
each $k$ we have
\begin{equation}
\lambda\left(Q(z,\rho_{k})\right)<\lambda(V^{\circ}\setminus F_{+})\label{eq:QVzero}
\end{equation}
 and

\begin{equation}
Q(z,\rho_{k})\subset V^{\circ}\label{eq:NearlyForg}
\end{equation}
and 
\[
\lambda\left(Q(z,\rho_{k})\cap F_{+}\right)<\lambda\left(Q(z,\rho_{k})\right)
\]
 and therefore also

\begin{equation}
\lambda\left(Q(z,\rho_{k})\setminus F_{+}\right)>0.\label{eq:mnky}
\end{equation}
for each $k\in\mathbb{N}$.

Let us now define a sequence $\left\{ F_{+}(k)\right\} _{k\in\mathbb{N}}$
of measurable sets by setting 
\begin{equation}
F_{+}(k)=F_{+}\cup Q(z,\rho_{k})\mbox{ for each\,}k\in\mathbb{N}.\label{eq:FplusK}
\end{equation}
 Note that, in view of \eqref{eq:mnky},
\begin{equation}
\lambda\left(F_{+}(k)\right)=\lambda(F_{+})+\lambda\left(Q(z,\rho_{k})\setminus F_{+}\right)>\lambda(F_{+}).\label{eq:Grermpst}
\end{equation}
We also introduce the set 
\begin{equation}
G:=F_{-}\setminus V^{\circ}.\label{eq:defGVFmin}
\end{equation}
 Note that 
\begin{equation}
F_{+}(k)\cap G=\emptyset\text{ for each\,}k\in\mathbb{N}\label{eq:FkGdisjoint}
\end{equation}
since in fact 
\begin{eqnarray*}
F_{+}(k)\cap G & = & \left(F_{+}\cup Q(z,\rho_{k})\right)\cap G=(F_{+}\cap G)\cup(Q(z,\rho_{k})\cap G)\\
 & \subset & (F_{+}\cap F_{-})\cup\left(Q(z,\rho_{k})\setminus V^{\circ}\right)=\emptyset\cup\emptyset.
\end{eqnarray*}

We also claim that 
\begin{equation}
\lambda\left(\left[0,1\right]^{d}\setminus F_{+}(k)\setminus G\right)>0\text{ for all\,}k\in\mathbb{N}.\label{eq:WWCN}
\end{equation}
To show this we first observe that, obviously, $V^{\circ}\setminus F_{+}(k)=V^{\circ}\setminus F_{+}(k)\setminus G$.
Consequently, 
\begin{eqnarray*}
 &  & \lambda\left(\left[0,1\right]^{d}\setminus F_{+}(k)\setminus G\right)\ge\lambda\left(V^{\circ}\setminus F_{+}(k)\setminus G\right)=\lambda\left(V^{\circ}\setminus F_{+}(k)\right)\\
 & = & \lambda\left(V^{\circ}\setminus F_{+}\setminus Q(z,\rho_{k})\right)=\lambda\left(V^{\circ}\setminus F_{+}\right)-\lambda\left((V^{\circ}\setminus F_{+})\cap Q(z,\rho_{k})\right)\\
 & \ge & \lambda\left(V^{\circ}\setminus F_{+}\right)-\lambda\left(Q(z,\rho_{k})\right).
\end{eqnarray*}
By \eqref{eq:QVzero} this last expression is strictly positive for
each $k$, which completes our proof of \eqref{eq:WWCN}.

We next observe that, by the first part of \eqref{eq:OmegaStuff}
(and of course also \eqref{eq:defGVFmin}),
\begin{equation}
\lambda(F_{-})=\lambda(F_{-}\cap V^{\circ})+\lambda(F_{-}\setminus V^{\circ})=\lambda(F_{-}\setminus V^{\circ})=\lambda(G).\label{eq:FGstuff}
\end{equation}

We use \eqref{eq:WWCN} then \eqref{eq:FkGdisjoint} then \eqref{eq:Grermpst}
and \eqref{eq:FGstuff}, and finally the disjointness of $F_{+}$
and $F_{-}$ to obtain that 
\begin{eqnarray*}
0 & < & \lambda\left([0,1]^{d}\setminus F_{+}(k)\setminus G\right)=1-\lambda(F_{+}(k))-\lambda(G)\\
 & < & 1-\lambda\left(F_{+}\right)-\lambda(F_{-})=\lambda\left(\left[0,1\right]^{d}\setminus F_{+}\setminus F_{-}\right).
\end{eqnarray*}
Since $\tau>0$, this implies that 
\begin{equation}
0<\tau\lambda\left([0,1]^{d}\setminus F_{+}(k)\setminus G\right)<\tau\lambda\left([0,1]^{d}\setminus F_{+}\setminus F_{-}\right).\label{eq:ilfd}
\end{equation}
Then, since $\min\left\{ \lambda(F_{+}(k)),\lambda(G)\right\} \ge\min\left\{ \lambda(F_{+}),\lambda(F_{-})\right\} $,
it follows from \ref{eq:ilfd} and the fact that $F_{+}$ and $F_{-}$
satisfy \eqref{eq:EqualityAllowed} that 
\[
\min\left\{ \lambda(F_{+}(k)),\lambda(G)\right\} >\tau\lambda(\left[0,1\right]^{d}\setminus F_{+}(k)\setminus G)>0\mbox{ for each\,}k\in\mathbb{N}.
\]
We use these strict inequalities to construct a sequence $\left\{ j(k)\right\} _{k\in\mathbb{N}}$
of positive integers such that $j(k)\ge k$ and $j(k)$ is sufficiently
large to ensure that $\tau_{j(k)}$ is sufficiently close to $\tau$
to imply that 
\[
\min\left\{ \lambda(F_{+}(k)),\lambda(G))\right\} >\tau_{j(k)}\lambda(\left[0,1\right]^{d}\setminus F_{+}(k)\setminus G)\mbox{ for each\,}k\in\mathbb{N}.
\]
Since $\left(\tau_{j(k)},s_{j(k)}\right)\in JS(d)$ and since $F_{+}(k)$
and $G$ are disjoint, we see, in accordance with Definition \ref{def:JohnStrombergPair}
(for $\mathcal{E}=\mathcal{Q}(\mathbb{R}^{d})$), that there exists
a subcube $W(k)$ of $\left[0,1\right]^{d}$ which satisfies 
\begin{equation}
\min\left\{ \lambda(W(k)\cap F_{+}(k)),\lambda(W(k)\cap G)\right\} \ge s_{j(k)}\lambda(W(k))\,.\label{eq:Chomp}
\end{equation}
The fact that $j(k)\ge k$ for each $k$ ensures that
\begin{equation}
\lim_{k\to\infty}s_{j(k)}=s.\label{eq:S-Limit}
\end{equation}

\noindent Since 
\[
W(k)\cap F_{+}\subset W(k)\cap F_{+}(k)=(W(k)\cap F_{+})\cup\left(W(k)\cap Q(z,r_{k})\right)\subset(W(k)\cap F_{+})\cup Q(z,r_{k})
\]

\noindent and $\lambda\left(Q(z,r_{k})\right)=\left(2r_{k}\right)^{d}$,
we obtain that 
\begin{equation}
\lambda\left(W(k)\cap F_{+}\right)\le\lambda\left(W(k)\cap F_{+}(k)\right)\le\lambda(W(k)\cap F_{+})+\left(2\rho_{k}\right)^{d}\mbox{ for all\,\ }k\in\mathbb{N}.\label{eq:SandwichA}
\end{equation}
A slight variant of the simple reasoning in \eqref{eq:FGstuff}, again
using the first part of \eqref{eq:OmegaStuff} and \eqref{eq:defGVFmin},
gives us that
\begin{align}
\lambda\left(W(k)\cap F_{-}\right) & =\lambda\left(W(k)\cap F_{-}\cap V^{\circ}\right)+\lambda\left(W(k)\cap F_{-}\setminus V^{\circ}\right)\nonumber \\
 & \phantom{=}\label{eq:FGZstufBIS}\\
 & =\lambda\left(W(k)\cap F_{-}\setminus V^{\circ}\right)=\lambda\left(W(k)\cap G\right)\text{\,for all\,}k\in\mathbb{N}.\nonumber 
\end{align}

For each $k\in\mathbb{N}$ we let the point $x_{k}$ in $\left[0,1\right]^{d}$
and the number $r_{k}\in(0,1/2]$ be the centre and half-side length
respectively of the subcube $W(k)$. I.e., we have $W(k)=Q(x_{k},r_{k})$.
There exists a strictly increasing sequence $\left\{ n_{k}\right\} _{k\in\mathbb{N}}$
of positive integers such that the sequences $\left\{ x_{n_{k}}\right\} _{k\in\mathbb{N}}$
and $\left\{ r_{n_{k}}\right\} _{k\in\mathbb{N}}$ converge, respectively,
to a point $x\in[0,1]^{d}$ and a number $r\in[0,1/2]$. We can now
declare the sequence $\left\{ W(n_{k})\right\} _{k\in\mathbb{N}}=\left\{ Q\left(x_{n_{k}},r_{n_{k}}\right)\right\} _{k\in\mathbb{N}}$
to be the special convergent sequence which we set out to construct
in this step of the proof in the case where $\left(F_{+},F_{-}\right)$
satisfies condition (ii) of Definition \ref{def:Tame}.

Since we have already specified our construction of $\left\{ W(n_{k})\right\} _{k\in\mathbb{N}}$
when $\left(F_{+},F_{-}\right)$ satisfies condition (i) of Definition
\ref{def:Tame}, this completes Step 1 of our proof.

\smallskip{}

\textbf{\textit{STEP 2: A proof that whenever the limit of $\left\{ W(n_{k})\right\} _{k\in\mathbb{N}}$
is a cube, then that cube has the two properties required to immediately
complete the proof of the theorem.}}

Suppose that the limit of the sequence $\left\{ W(n_{k})\right\} _{k\in\mathbb{N}}=\left\{ Q\left(x_{n_{k}},r_{n_{k}}\right)\right\} _{k\in\mathbb{N}}$
which was constructed in the previous step, is indeed a cube $W=Q(x,r)$,
i.e., that 
\begin{equation}
r:=\lim_{k\to\infty}r_{n_{k}}>0.\label{eq:rpos}
\end{equation}
Let us now prove that this implies that $W$ is contained in $\left[0,1\right]^{d}$
and that it satisfies \eqref{eq:Bloomp}. (We defer treatment of the
case where $\lim_{k\to\infty}r_{n_{k}}=0$ to Step 4.)

In the case which was dealt with at the beginning of the previous
step, where $F_{+}$ and $F_{-}$ satisfy condition (i) of Definition
\ref{def:Tame}, and where indeed we always have $r>0$, it is already
known that the cube $W$ is contained in $\left[0,1\right]^{d}.$
We note that, by the reasoning in Remark \ref{rem:s-half}, we have
$s_{n}\le1/2$ for each $n\in\mathbb{N}$ and therefore also $s\le1/2$.
So \eqref{eq:Bloomp} follows immediately from \eqref{eq:EasyStuff}
and \eqref{eq:same}.

We turn to the remaining case, where $F_{+}$ and $F_{-}$ satisfy
condition (ii) of Definition \ref{def:Tame}, and therefore the sequence
$\left\{ W\left(n_{k}\right)\right\} _{k\in\mathbb{N}}$ is constructed
in the more elaborate way described in the second and much longer
part of the previous step. Here the positivity of the limit $r$ permits
us to apply Lemma 10.1 of \cite[p.~43]{CwikelMSagherYShvartsmanP2010}
(whose proof was briefly discussed above near the end of the first
paragraph of the proof of Theorem \ref{thm:onlyneedcubes}) to the
sequence $\left\{ W(n_{k})\right\} _{k\in\mathbb{N}}=\left\{ Q(x_{n_{k}},r_{n_{k}})\right\} _{k\in\mathbb{N}}$
to obtain the first required conclusion, that $W\subset\left[0,1\right]^{d}$,
and to also obtain that 
\begin{equation}
\lambda(W\cap F_{+})=\lim_{k\to\infty}\lambda\left(W(n_{k})\cap F_{+}\right)\label{eq:ruaz}
\end{equation}
and (here also using \eqref{eq:FGZstufBIS}) that 
\begin{equation}
\lambda(W\cap F_{-})=\lim_{k\to\infty}\lambda\left(W(n_{k})\cap F_{-}\right)=\lim_{k\to\infty}\lambda\left(W(n_{k})\cap G\right)=\lambda(W\cap G),\label{eq:ruaz-aux}
\end{equation}
and also that 
\begin{equation}
\lambda(W)=\lim_{k\to\infty}\lambda\left(W(n_{k})\right).\label{eq:Preimp}
\end{equation}
From \eqref{eq:ruaz} and \eqref{eq:SandwichA} we see that 
\begin{equation}
\lambda(W\cap F_{+})=\lim_{k\to\infty}\lambda\left(W(n_{k})\cap F_{+}(n_{k})\right),\label{eq:ging}
\end{equation}
and from \eqref{eq:Preimp} and \eqref{eq:S-Limit} we see that 
\begin{equation}
s\lambda(W)=\lim_{k\to\infty}s_{j(n_{k})}\lambda\left(W(n_{k})\right).\label{eq:nong}
\end{equation}
In view of \eqref{eq:Chomp} we have that 
\[
\ensuremath{\ensuremath{\lambda(W(n_{k})\cap F_{+}(n_{k}))\ge s_{j(n_{k})}\lambda(W(n_{k}))\mbox{ \,\ and\,\,}\lambda(W(n_{k})\cap G)\ge s_{j(n_{k})}\lambda(W(n_{k}))}}
\]
for each $k\in\mathbb{N}$. If we take the limit as $k$ tends to
$\infty$ in each of these two inequalities, and apply \eqref{eq:ging},
\eqref{eq:ruaz-aux} and \eqref{eq:nong}, then the two resulting
inequalities can be rewritten as the required single inequality \eqref{eq:Bloomp}.

Thus we have shown that if \eqref{eq:rpos} holds, then in both cases,
i.e., whether it is condition (i) or condition (ii) of Definition
\ref{def:Tame} which applies to $F_{+}$ and $F_{-}$, it follows
that the limiting cube $W$ indeed has both the properties required
to complete the proof of part (a) and also part (b) of the lemma.

It remains to explain, as we shall do in the remaining two steps of
the proof, how we can sometimes guarantee that \eqref{eq:rpos} does
hold, and how we can proceed when it does not hold.

\smallskip{}

\textbf{\textit{STEP 3: Completion of the proof of part (a) of the
lemma.}}

Our preceding treatment of the case where \eqref{eq:rpos} holds,
will now enable us to complete the proof of part (a) of the lemma
in full generality. Part (a) refers to the property $(**)$, so in
our proof of it we can and must assume that the disjoint sets $F_{+}$
and $F_{-}$ are both $d$-multi-cubes.

We once again refer (as we did at the beginning of the proof of Theorem
\ref{thm:EqualityAllowed}) to the simple reasoning after \eqref{eq:wwnt}
in the proof of Theorem \ref{thm:SmallerTau} which shows that there
exists a dyadic cube $V$ whose interior is contained in $\left[0,1\right]^{d}\setminus F_{+}\setminus F_{-}$
and which therefore must satisfy \eqref{eq:OmegaStuff}. So we can
construct the sequence of cubes $\left\{ W(n_{k})\right\} _{k\in\mathbb{N}}$
for this particular choice of $V$ in exactly the way that was done
in the second part of Step 1 of this proof. This choice of $V$ implies
that the set $G$ introduced in \eqref{eq:defGVFmin} must satisfy
$G=F_{-}\setminus V^{\circ}=F_{-}$.

Since they are $d$-multi-cubes, $F_{+}$ and $F_{-}=G$ are compact
and so of course is $Q(z,\rho_{k})$. Consequently, the set $F_{+}(k)$
(defined by \eqref{eq:FplusK}) is also compact. In view of \eqref{eq:NearlyForg},
it is disjoint from $G$. So the distance, which we denote by $\mathrm{dist_{\infty}}\left(F_{+}(k),G\right)$,
between $F_{+}(k)$ and $G$ with respect to the $\ell^{\infty}$
metric on $\mathbb{R}^{d}$, must be positive. In fact, since $F_{+}(k)\subset F_{+}(1)$
(because $\rho_{k}\le\rho_{1})$, we have that 
\begin{equation}
\mathrm{dist_{\infty}}\left(F_{+}(k),G\right)\ge\mathrm{dist_{\infty}}\left(F_{+}(1),G\right)>0\,\,\,\mbox{for all\,}k\in\mathbb{N}.\label{eq:PosDist}
\end{equation}

Now we proceed more or less similarly to the last steps of the proof
of Theorem \ref{thm:SmallerTau}. By \eqref{eq:Chomp} the cube $W(k)=Q(x_{k},r_{k})$
intersects both of the sets $F_{+}(k)$ and $G$. So its side length
$2r_{k}$ cannot be smaller than $\mathrm{dist_{\infty}}\left(F_{+}(k),G)\right)$.
Consequently, also using \eqref{eq:PosDist}, we see that 
\[
\frac{1}{2}\ge r_{k}\ge\delta_{0}:=\frac{1}{2}\mathrm{dist_{\infty}}\left(F_{+}(1),G\right)>0\mbox{ for all\,}k\in\mathbb{N}.
\]
Consequently $r\ge\delta_{0}$ and so \eqref{eq:rpos} holds, permitting
us to use the reasoning of Step 2 to ensure the existence of cube
$W$ which has the properties required to complete the proof of part
(a).

\smallskip{}

\textbf{\textit{STEP 4: Completion of the proof of part (b) of the
lemma.}}

In view of Step 2, if the sequence of cubes $\left\{ W(n_{k})\right\} _{k\in\mathbb{N}}$
constructed in Step 1 converges to the cube $W$, then that cube is
contained in $\left[0,1\right]^{d}$ and satisfies \eqref{eq:Bloomp}
and no further reasoning is required to complete the proof of part
(b). Thus it remains only to deal with the case where $r=\lim_{k\to\infty}r_{n_{k}}=0$.
It is clear from the first part of Step 1 that this cannot happen
if $\left(F_{+},F_{-}\right)$ satisfies condition (i) of Definition
\ref{def:Tame}. So the sequence $\left\{ W(n_{k})\right\} _{k\in\mathbb{N}}$
has necessarily been constructed as in the second part of Step 1,
via the sets $F(k)$ and other sets introduced there.

Let us first see that, in this case, the point $x=\lim_{k\to\infty}x_{n_{k}}$
cannot coincide with the point $z$ which appears in the definition
\eqref{eq:FplusK} of the sets $F_{+}(k)$. If $x=z$ and is therefore
in the interior $V^{\circ}$ of the dyadic cube $V$ used in the construction,
then there exists some $k$ which is sufficiently large to ensure
that $r_{n_{k}}+\left\Vert x_{n_{k}}-x\right\Vert _{\infty}<\mathrm{dist}_{\infty}\left(x,\partial V\right)$.
This means that every point $y$ in the cube $W(n_{k})=Q(x_{n_{k}},r_{n_{k}})$
satisfies $\left\Vert y-x\right\Vert _{\infty}\le\left\Vert y-x_{n_{k}}\right\Vert _{\infty}+\left\Vert x_{n_{k}}-x\right\Vert _{\infty}<\mathrm{dist}_{\infty}\left(x,\partial V\right)$.
Consequently $W(n_{k})\subset V^{\circ}$ and therefore $W(n_{k})\cap G=\emptyset$.
(Recall that $G$ is defined by \eqref{eq:defGVFmin}). Since this
contradicts \eqref{eq:Chomp}, we have indeed shown that $x\ne z$.

In view of this fact, there exists an integer $k_{0}$ which is large
enough to ensure that 
\begin{equation}
\left\Vert x-x_{n_{k_{0}}}\right\Vert _{\infty}+\rho_{n_{k_{0}}}+r_{n_{k_{0}}}<\left\Vert x-z\right\Vert _{\infty},\label{eq:rhor}
\end{equation}
where $\rho_{n_{k_{0}}}$ is a element of the sequence $\left\{ \rho_{k}\right\} _{k\in\mathbb{N}}$
with limit $0$ which is used in the definition \eqref{eq:FplusK}
of the sets $F_{+}(k)$. If $y$ is a point in the intersection of
the two cubes $W(n_{k_{0}})=Q\left(x_{n_{k_{0}}},r_{n_{k_{0}}}\right)$
and $Q\left(z,\rho_{n_{k_{0}}}\right)$, then 
\begin{eqnarray*}
\left\Vert x-z\right\Vert _{\infty} & \le & \left\Vert x-x_{n_{k_{0}}}\right\Vert _{\infty}+\left\Vert x_{n_{k_{0}}}-y\right\Vert _{\infty}+\left\Vert y-z\right\Vert _{\infty}\\
 & \le & \left\Vert x-x_{n_{k_{0}}}\right\Vert _{\infty}+\rho_{n_{k_{0}}}+r_{n_{k_{0}}}.
\end{eqnarray*}
But this contradicts \eqref{eq:rhor} and enables us to conclude that
$W\left(n_{k_{0}}\right)$ and $Q\left(z,\rho_{n_{k_{0}}}\right)$
must be disjoint and therefore that 
\[
W\left(n_{k_{0}}\right)\cap F_{+}\left(n_{k_{0}}\right)=W\left(n_{k_{0}}\right)\cap F_{+}.
\]
We apply this, together with \eqref{eq:FGZstufBIS} for $k=n_{k_{0}}$
and then \eqref{eq:Chomp} for $k=n_{k_{0}}$, to obtain that the
cube $W:=W\left(n_{k_{0}}\right)$ satisfies 
\begin{eqnarray}
\min\left\{ \lambda\left(W\cap F_{+}\right),\lambda\left(W\cap F_{-}\right)\right\}  & = & \min\left\{ \lambda\left(W\cap F_{+}\left(n_{k_{0}}\right)\right),\lambda\left(W\cap G\right)\right\} \nonumber \\
 & \ge & s_{j(n_{k_{0}})}\lambda(W).\label{eq:oomze}
\end{eqnarray}
The cube $W$, like all other cubes in the sequence $\left\{ W(n_{k})\right\} _{k\in\mathbb{N}}$
is contained in $\left[0,1\right]^{d}$. Finally we have to recall
that in the statement of part (b) of the lemma, the sequence $\left\{ s_{n}\right\} _{n\in\mathbb{N}}$
is required to satisfy $s_{n}\ge s$ for each $n$. So \eqref{eq:oomze}
shows that, also in this last remaining case, we have obtained a subcube
$W$ of $\left[0,1\right]^{d}$ which satisfies \eqref{eq:Bloomp}
for the given sets $F_{+}$ and $F_{-}$. This completes the proof
of the lemma. $\qed$

\subsection{Some pairs which are not John-Str\"omberg pairs}

The following result is a more elaborate variant of Remark \ref{rem:s-half}.
\begin{lem}
\label{lem:NotInJSd}For each $d\in\mathbb{N}$ and each $\tau>0$,
and for every $s>1/(2+1/\tau)$, the pair $\left(\tau,s\right)$ is
\textbf{not} in $JS(d)$.
\end{lem}

\noindent \textit{Proof.} If $s>1/2$ then the result follows from
Remark \ref{rem:s-half}. So we can assume that $1/(2+1/\tau)<s\le1/2$.
Let us choose some number $a\in\left(1/(2+1/\tau),s\right)$ and then
let $E_{-}=\left\{ (x,t):x\in[0,1]^{d-1},0\le t\le a\right\} $ and
$E_{+}=\left\{ (x,t):x\in[0,1]^{d-1},1-a\le t\le1\right\} $. (Obviously
when $d=1$ we have to interpret the previous definition to mean that
$E_{-}=[0,a]$ and $E_{+}=[1-a,1].)$ Since $a<1/2$, these two sets
are disjoint. Furthermore 
\[
\lambda\left(\left[0,1\right]^{d}\setminus E_{+}\setminus E_{-}\right)=1-2a=\frac{1-2a}{a}\min\left\{ \lambda(E_{+}),\lambda\left(E_{-}\right)\right\} 
\]
or, equivalently, 
\[
\min\left\{ \lambda(E_{+}),\lambda\left(E_{-}\right)\right\} =\frac{1}{\frac{1}{a}-2}\lambda\left(\left[0,1\right]^{d}\setminus E_{+}\setminus E_{-}\right).
\]
Since $a>1/(2+1/\tau)$ it follows that $\frac{1}{\frac{1}{a}-2}>\tau$
and so 
\begin{equation}
\min\left\{ \lambda(E_{+}),\lambda\left(E_{-}\right)\right\} >\tau\lambda\left(\left[0,1\right]^{d}\setminus E_{+}\setminus E_{-}\right).\label{eq:gpyr}
\end{equation}
We will complete the proof of this lemma by showing that, although
the two disjoint measurable subsets $E_{+}$ and $E_{-}$ of $\left[0,1\right]^{d}$
satisfy \eqref{eq:gpyr}, there does not exist any subcube $W$ of
$\left[0,1\right]^{d}$ which satisfies $\min\left\{ \lambda(E_{+}\cap W),\lambda\left(E_{-}\cap W\right)\right\} \ge s$$\lambda(W)$.
This will follow from the inequality 
\begin{equation}
\frac{\min\left\{ \lambda(E_{+}\cap W),\lambda\left(E_{-}\cap W\right)\right\} }{\lambda(W)}\le a\label{eq:ermp}
\end{equation}
which will be seen to hold for every subcube $W$ of $\left[0,1\right]^{d}$.
This inequality seems intuitively quite obvious, but let us nevertheless
give a detailed (and perhaps not optimally elegant) proof.

Let $W$ be an arbitrary subcube of $\left[0,1\right]^{d}$ and let
$\theta$ be its side length. Of course $\theta\in(0,1]$ and $W$
must be of the form $W=\left\{ (x,t):x\in W_{0},t\in[\beta,\beta+\theta]\right\} $
where $W_{0}$ is some subcube of $\left[0,1\right]^{d-1}$ of side
length $\theta$ and $\left[\beta,\beta+\theta\right]\subset\left[0,1\right]$.
(If $d=1$ then $W$ is simply the interval $\left[\beta,\beta+\theta\right]$.)
We only need to consider the case where $\left[\beta,\beta+\theta\right]$
has a non-empty intersection with each one of the intervals $\left[0,a\right]$
and $\left[1-a,1\right]$, since otherwise at least one of the two
sets $E_{+}\cap W$ and $E_{-}\cap W$ is empty and \eqref{eq:ermp}
is a triviality. The non-emptiness of the above mentioned two intersections
implies that 
\[
\beta\le a\mbox{ and\,}\beta+\theta\ge1-a.
\]
When $\theta$ and $\beta$ satisfy all the above mentioned conditions
we have that $\lambda(E_{-}\cap W)=\theta^{d-1}(a-\beta)$ and $\lambda\left(E_{+}\cap W\right)=\theta^{d-1}\left(\beta+\theta-(1-a)\right)$.
Therefore 
\begin{equation}
\min\left\{ \lambda(E_{+}\cap W),\lambda\left(E_{-}\cap W\right)\right\} \le\theta^{d-1}\sup_{\beta\in\mathbb{R}}\left[\min\left\{ a-\beta,\beta+\theta-(1-a)\right\} \right].\label{eq:hqr}
\end{equation}
For each fixed choice of $\theta$, the expression in the square brackets
on the right side of \eqref{eq:hqr} is the minimum of a strictly
decreasing function of $\beta$ and a strictly increasing function
of $\beta$. Therefore its supremum and thus its maximum is attained
when $\beta$ takes the unique value for which these two functions
are equal, namely when $\beta=\left(1-\theta\right)/2$. This shows
that 
\begin{equation}
\frac{\min\left\{ \lambda(E_{+}\cap W),\lambda\left(E_{-}\cap W\right)\right\} }{\lambda(W)}\le\frac{\theta^{d-1}}{\theta^{d}}\left(a-\frac{1-\theta}{2}\right)=\frac{1}{2}-\frac{1}{\theta}\left(\frac{1}{2}-a\right).\label{eq:hong}
\end{equation}
Finally, the facts that $\theta\in(0,1]$ and $\frac{1}{2}-a>0$ imply
that the right side of \eqref{eq:hong} is bounded above by 
\[
\frac{1}{2}-\left(\frac{1}{2}-a\right)=a
\]
which establishes \eqref{eq:ermp} and so completes the proof of the
lemma. $\qed$
\begin{rem}
\label{rem:RonIsOptimal} In particular, when $\tau=1/2$, Lemma \ref{lem:NotInJSd}
shows that 
\begin{equation}
\left(\frac{1}{2},\frac{1}{4}+\varepsilon\right)\notin JS(d)\text{\,\,for every\,\,}d\in\mathbb{N\text{\,\,and every\,\,}}\varepsilon>0.\label{eq:BestPossible}
\end{equation}
Thus the second of the three results obtained recently by Ron Holzman
(see Section \ref{sec:KnownResults}), namely that $\left(\frac{1}{2},\frac{1}{4}\right)\in JS(1)$,
shows that $s=\frac{1}{4}$ is the largest possible value of $s$
for which $\left(\frac{1}{2},s\right)\in JS(1)$. In the light of
\eqref{eq:BestPossible}, the first of his three results can also
be considered, in some sense, to be best possible. At least in the
non-trivial special case that he considered there he obtained that
there exists a cube $W$ in $\left[0,1\right]^{d}$ which satisfies
\eqref{eq:blurp} for $s=1/4$. (Apparently one cannot exclude the
possibility that, at least in that special case, one might be able
to also obtain a cube $W$ in $\left[0,1\right]^{d}$ which satisfies
\eqref{eq:blurp} for some $s>1/4$.) These results tempt one to wonder
whether perhaps the property $\left(\frac{1}{2},\frac{1}{4}\right)\in JS(d)$
might hold for \textit{all} $d\in\mathbb{N}$, If so, that would of
course also answer Question A(1/2) affirmatively, and in an optimally
strong, even ``dramatically strong'' way. However, in the third
of his results relating to this question, Holtzman has analysed the
following example suggested by the author and has shown that this
is an ``impossible dream''. This property fails to hold already
for $d=2$. Therefore (cf.~Theorem \ref{thm:DplusOne} below) it
also does not hold for any other $d>1$.

Let $F_{+}$ be the rectangle $F_{+}=$$\left[0,1\right]\times\left[2/3,1\right]$
and let $F_{-}$ be the union of the two squares $\left[0,1/3\right]\times\left[0,1/3\right]$
and $\left[2/3,1\right]\times\left[0,1/3\right]$. Then $\left(F_{+},F_{-}\right)$
is a tame couple in $\left[0,1\right]^{2}$. The areas of the three
disjoint sets $F_{+},$ $F_{-}$ and $\left[0,1\right]^{2}\setminus F_{+}\setminus F_{-}$
are respectively $1/3$, 2/9 and $4/9$, and this ensures that $F_{+}$
and $F_{-}$ satisfy \eqref{eq:EqualityAllowed} for $\tau=1/2$ and
$d=2$. (In this example $\lambda$ will of course always denote two-dimensional
Lebesgue measure.) Let $\mathcal{W}$ denote the collection of all
closed squares with sides parallel to the axes which are contained
in $\left[0.1\right]^{2}$. For each $W\in\mathcal{W}$ let 
\[
f(W)=\frac{\min\left\{ \lambda(W\cap F_{+}),\lambda(W\cap F_{-})\right\} }{\lambda(W)}.
\]
The third of Ron Holzman's results is that, for these choices of $F_{+}$
and $F_{-}$, 
\[
\sup\left\{ f(W):W\in\mathcal{W}\right\} =\max\left\{ f(W):W\in\mathcal{W}\right\} =\sqrt{5}-2.
\]
In view of part $(***)$ of Theorem \ref{thm:EqualityAllowed}, this
shows that, indeed, $\left(1/2,1/4\right)\notin JS(2)$ and, furthermore
(again recalling Theorem \ref{thm:DplusOne}), that $\left(1/2,s\right)\notin JS(d)$
for every $s>\sqrt{5}-2$ and $d\ge2$.

It is tempting to wonder whether a sequence of appropriate variants
of this example for subsets of $\left[0,1\right]^{d}$, where $d$
tends to $\infty$, might lead to a negative answer to Question A(1/2).
Initial attempts to find such a sequence have not yielded anything
decisive.
\end{rem}

\subsection{Some further properties of the set $JS(d)$}
\begin{fact}
For each $d\in\mathbb{N}$ and each $\tau>0$ there exists some $s>0$
such that $\left(\tau,s\right)\in JS(d).$
\end{fact}

\noindent \textit{Proof.} This is a consequence of Theorem 7.8 of
\cite[pp.~28--29]{CwikelMSagherYShvartsmanP2010} a.k.a. Theorem 7.7
of \cite[pp.~152--153]{CwikelMSagherYShvartsmanP2012} which, for
each $\tau>0$, provides a positive number $s$ depending on $\tau$
and $d$ such that $\left(\tau,s\right)\in JS(d)$. (The explicit
formula for $s$ will be recalled and used in Theorem \ref{thm:PropsJSd}
below. $\qed$

The preceding result ensures that the supremum in the following definition
is taken over a non-empty set.
\begin{defn}
For each $d\in\mathbb{N}$ and $\tau>0$ let
\[
\sigma(\tau,d):=\sup\left\{ s>0:\left(\tau,s\right)\in JS(d)\right\} .
\]
\end{defn}

\begin{rem}
\label{rem:RonAgain}Thus the second and third results of Ron Holzman
mentioned in Section \ref{sec:KnownResults} (cf.~also Remark \ref{rem:RonIsOptimal})
can be written as
\[
\sigma\left(\frac{1}{2},1\right)=\frac{1}{4}\mbox{ \,\,\ and\,\,\,\,}\sigma\left(\frac{1}{2},2\right)\le\sqrt{5}-2.
\]
We can now readily establish several properties of the set $JS(d)$
and the function $\sigma(\tau,d)$.
\end{rem}

\begin{thm}
\label{thm:PropsJSd}For each fixed $d\in\mathbb{N},$

$\mathrm{(i)}$~$(\tau,\sigma(\tau,d)]\in JS(d)$ for each $\tau>0$.

$\mathrm{(ii)}$~ $JS(d)=\left\{ (x,y):x>0,\,0<y\le\sigma(x,d)\right\} $

$\mathrm{(iii)}$~ The function $x\mapsto\sigma(x,d)$ is non-decreasing
and continuous and satisfies 
\begin{equation}
\varphi(x,d)\le\sigma(x,d)\le\frac{1}{2+\frac{1}{x}}\mbox{ for all\,}x>0,\label{eq:Foomp}
\end{equation}
where 
\begin{equation}
\varphi(x,d)=\left.\left\{ \begin{array}{ccc}
2^{-d}(x-x^{2})/(1+x) & , & 0<x\le\sqrt{2}-1\\
2^{-d}\left(3-2\sqrt{2}\right) & , & x\ge\sqrt{2}-1.
\end{array}\right.\right.\label{eq:Form4Phi}
\end{equation}
\end{thm}

\noindent \textit{Proof.} ~For part (i) let us fix an arbitrary
$\tau>0$ and first note the obvious fact that 
\begin{equation}
\left(\tau,s\right)\in JS(d)\Rightarrow\left(\tau,s'\right)\in JS(d)\mbox{ for all\,}s'\in(0,s),\label{eq:ObvImp}
\end{equation}
which implies that $\left(\tau,s'\right)\in JS(d)$ for every $s'\in\left(0,\sigma(\tau,d)\right).$
Therefore the sequences $\left\{ \tau_{n}\right\} _{n\in\mathbb{N}}$
and $\left\{ s_{n}\right\} _{n\in\mathbb{N}}$ which we define by
$\tau_{n}:=\tau$ and $s_{n}:=\left(1-2^{-n}\right)\sigma(\tau,d)$
must satisfy $\left(\tau_{n},s_{n}\right)\in JS(d)$ for every $n\in\mathbb{N}$.
For the proof of part (i) we now simply apply part (a) of Lemma \ref{lem:Sequences of JS pairs}
to the sequences $\left\{ \tau_{n}\right\} _{n\in\mathbb{N}}$ and
$\left\{ s_{n}\right\} _{n\in\mathbb{N}}$ and then apply Theorem
\ref{thm:EqualityAllowed}. Part (ii) then follows immediately from
part (i) and \eqref{eq:ObvImp}. For part (iii) we first use another
obvious fact, namely that 
\[
\left(\tau,s\right)\in JS(d)\Rightarrow\left(\tau',s\right)\in JS(d)\mbox{ for all\,}\tau'>\tau
\]
to immediately deduce that the function $x\mapsto\sigma(x,d)$ is
non-decreasing. This latter property means that, in order to show
that this function is continuous, it will suffice to show that 
\begin{equation}
\lim_{n\to\infty}\sigma\left((1+1/n)\tau,d\right)\le\sigma(\tau,d)\mbox{ and\,}\lim_{n\to\infty}\sigma\left((1-1/n)\tau,d\right)\ge\sigma(\tau,d)\mbox{ for each\,}\tau>0.\label{eq:ForCty}
\end{equation}
We obtain the first of these inequalities by again using Theorem \ref{thm:EqualityAllowed}
together with part (a) of Lemma \ref{lem:Sequences of JS pairs} to
show that the limit of the sequence $\left\{ \left((1+1/n)\tau,\sigma\left((1+1/n)\tau,d\right)\right)\right\} _{n\in\mathbb{N}}$
of points in $JS(d)$ must also be a point in $JS(d)$. We next remark
that, since $\left(\tau,\sigma(\tau,d)\right)\in JS(d)$, it follows
from Theorem \ref{thm:SmallerTau}, that $\left((1-1/n)\tau,\sigma(\tau,d)-\left(1-\sigma(\tau,d)\right)/n\right)\in JS(d)$
for all sufficiently large $n\in\mathbb{N}$. Therefore $\sigma\left((1-1/n)\tau,d\right)\ge\sigma(\tau,d)-\left(1-\sigma(\tau,d)\right)/n$
for these same values of $n$. This suffices to prove the second inequality
in \eqref{eq:ForCty} and complete the proof of continuity.

The formula \eqref{eq:Form4Phi} for $\varphi(x,d)$ in the estimate
from below in \eqref{eq:Foomp} is obtained by once more appealing
to Theorem 7.8 of \cite[pp.~28--29]{CwikelMSagherYShvartsmanP2010}
a.k.a. Theorem 7.7 of \cite[pp.~152--153]{CwikelMSagherYShvartsmanP2012},
and using the formulae appearing at the end of the statement of that
theorem in the particular case where the collection of sets $\mathcal{E}$
is chosen to be $\mathcal{Q}(\mathbb{R}^{d})$ the collection of all
cubes in $\mathbb{R}^{d}$. We can replace $M$ in those formulae
by $2^{d}$ using the fact that, in the terminology introduced just
before the statement of that theorem, $\mathcal{Q}(\mathbb{R}^{d})$
is $2^{d}$-decomposable. We can also replace $\delta$ there by $1/2$
using the fact (see Definition 7.4 of \cite[p.~26]{CwikelMSagherYShvartsmanP2010}
or \cite[p.~151]{CwikelMSagherYShvartsmanP2012} and the remark immediately
following it) that $1/2$ is a bi-density constant for $\mathcal{Q}(\mathbb{R}^{d})$.
The estimate from above in \eqref{eq:Foomp} follows from Lemma \ref{lem:NotInJSd}.

This completes the proof of part (iii) and therefore of the whole
theorem. $\qed$
\begin{rem}
There is another different kind of lower bound for $\sigma(\tau,d)$,
in terms of values of $\sigma(\tau',d)$ for appropriate numbers $\tau'$
greater than $\tau$, which can be obtained from Theorem \ref{thm:SmallerTau}.
We have not bothered to explicitly state it here.
\end{rem}

The following result seems intuitively completely obvious. But we
shall provide an explicit proof.
\begin{thm}
\label{thm:DplusOne}The inclusion $JS(d+1)\subset JS(d)$ and consequently
also the inequality $\sigma(\tau,d+1)\le\sigma(\tau,d)$ both hold
for every $d\in\mathbb{N}$ and $\tau>0$.
\end{thm}

\noindent \textit{Proof.} This is the one place in this paper where
we need to use the more explicit notation $\lambda_{d}$ instead of
$\lambda$ to denote $d$-dimensional Lebesgue measure. We will use
Theorem \ref{thm:onlyneedcubes}.

Suppose that $\left(\tau,s\right)\in JS(d+1)$. Let $F_{+}$ and $F_{-}$
be two \textit{arbitrarily chosen} disjoint $d$-multi-cubes which
satisfy \eqref{eq:Strict}, i.e., the inequality 
\begin{equation}
\min\left\{ \lambda_{d}(F_{+}),\lambda_{d}(F_{-})\right\} >\tau\lambda_{d}(\left[0,1\right]^{d}\setminus F_{+}\setminus F_{-})\,,\label{eq:vermeer}
\end{equation}
for this given value of $\tau$. We define two subsets $H_{+}$ and
$H_{-}$ of $\left[0,1\right]^{d+1}$ as the cartesian products $H_{+}=F_{+}\times[0,1]$
and $H_{-}=F_{-}\times[0,1]$. They are disjoint, since $F_{+}$ and
$F_{-}$are disjoint. For each dyadic subcube $E$ of $[0,1]^{d}$,
the cartesian product $E\times[0,1]$ is the union of $2^{k}$ dyadic
subcubes of $\left[0,1\right]^{d+1}$, where $k$ is such that the
side length of $E$ is $2^{-k}$. It follows that $H_{+}$ and $H_{-}$
are both $\left(d+1\right)$-multi-cubes.

Now, and also again later, we shall use the very standard facts that

(i) for each bounded closed interval $\left[a,b\right]$, the set
$E\times[a,b]$ is a Lebesgue measurable subset of $\mathbb{R}^{d+1}$
whenever $E$ is a Lebesgue measurable subset of $\mathbb{R}^{d}$,
and that

(ii)~~$\lambda_{d+1}(E\times[a,b])=(b-a)\lambda_{d}(E)$ for each
such $E$.

These two facts together with the fact that $\left[0,1\right]^{d+1}\setminus H_{+}\setminus H_{-}=(\left[0,1\right]^{d}\setminus F_{+}\setminus F_{-})\times[0,1]$,
and together with our assumption that $F_{+}$ and $F_{-}$ satisfy
\eqref{eq:vermeer}, imply that 
\[
\min\left\{ \lambda_{d+1}(H_{+}),\lambda_{d+1}(H_{-})\right\} >\tau\lambda_{d+1}(\left[0,1\right]^{d+1}\setminus H_{+}\setminus H_{-}).
\]
 Consequently, by Theorem \ref{thm:onlyneedcubes}, our assumption
that $\left(\tau,s\right)\in JS(d+1)$ guarantees the existence of
a subcube $V$ of $\left[0,1\right]^{d+1}$ for which 
\begin{equation}
\min\left\{ \lambda_{d+1}(V\cap H_{+}),\lambda_{d+1}(V\cap H_{-})\right\} \ge s\lambda_{d+1}(V).\label{eq:Vcube}
\end{equation}
We can write $V$ as the cartesian product $V=W\times[a,b]$, where
$W$ is a subcube of $\left[0,1\right]^{d}$ and $\left[a,b\right]$
is a subinterval of $\left[0,1\right]$ whose length $b-a$ of course
equals the side length of $V$ and of $W$. Obviously $V\cap H_{+}=\left(W\cap F_{+}\right)\times[a,b]$
and $V\cap H_{-}=\left(W\cap F_{-}\right)\times[a,b]$. So we can
again apply the standard facts (i) and (ii), which were recalled in
an earlier step of this proof, to obtain the formulae $\lambda_{d+1}(E)=(b-a)\lambda_{d}(E)$
in the three cases where $E$ is $W\cap F_{+}$ or $W\cap F_{-}$
or $W$. When we substitute these formulae in \eqref{eq:Vcube} and
divide both sides of the inequality by $b-a$, we obtain that the
subcube $W$ of $\left[0,1\right]^{d}$ satisfies \eqref{eq:nezt}
of Theorem \ref{thm:onlyneedcubes}, i.e., that
\[
\min\left\{ \lambda_{d}(W\cap F_{+}),\lambda_{d}(W\cap F_{-})\right\} \ge s\lambda_{d}(W)\,.
\]
Since $F_{+}$ and $F_{-}$ were chosen arbitrarily, we can once more
apply Theorem \ref{thm:onlyneedcubes} to deduce that $\left(\tau,s\right)\in JS(d),$
and so complete the proof of the present theorem. $\qed$

\section{\label{sec:CSS-Comments}Some comments and minor corrections for
the papers \cite{CwikelMSagherYShvartsmanP2010,CwikelMSagherYShvartsmanP2012}.}

In the first paragraph of the proof of Theorem 7.8 \cite[p.~30]{CwikelMSagherYShvartsmanP2010}
a.k.a. Theorem 7.7 \cite[pp.~154--155]{CwikelMSagherYShvartsmanP2012}
it is shown that it suffices to consider the case where the two sets
$E_{+}$ and $E_{-}$ are both compact. The justification of this
is a little clearer if in the formula on the third line of the proof
one replaces $G$ by $Q\setminus H_{+}\setminus H_{-},$ which is
obviously permissible.

We refer to the remarks made above in the course of the proof of Theorem
\ref{thm:onlyneedcubes} for some other small corrections and clarifications
of some small issues in \cite{CwikelMSagherYShvartsmanP2010}.

\end{document}